\documentclass[12pt]{article}
\usepackage{latexsym,amssymb,amsmath}
\begin{document}
\baselineskip=20pt

\newcommand{\la}{\langle}
\newcommand{\ra}{\rangle}
\newcommand{\psp}{\vspace{0.4cm}}
\newcommand{\pse}{\vspace{0.2cm}}
\newcommand{\ptl}{\partial}
\newcommand{\dlt}{\delta}
\newcommand{\sgm}{\sigma}
\newcommand{\al}{\alpha}
\newcommand{\be}{\beta}
\newcommand{\G}{\Gamma}
\newcommand{\gm}{\gamma}
\newcommand{\vs}{\varsigma}
\newcommand{\Lmd}{\Lambda}
\newcommand{\lmd}{\lambda}
\newcommand{\td}{\tilde}
\newcommand{\vf}{\varphi}
\newcommand{\yt}{Y^{\nu}}
\newcommand{\wt}{\mbox{wt}\:}
\newcommand{\rd}{\mbox{Res}}
\newcommand{\ad}{\mbox{ad}}
\newcommand{\stl}{\stackrel}
\newcommand{\ol}{\overline}
\newcommand{\ul}{\underline}
\newcommand{\es}{\epsilon}
\newcommand{\dmd}{\diamond}
\newcommand{\clt}{\clubsuit}
\newcommand{\vt}{\vartheta}
\newcommand{\ves}{\varepsilon}
\newcommand{\dg}{\dagger}
\newcommand{\tr}{\mbox{Tr}}
\newcommand{\ga}{{\cal G}({\cal A})}
\newcommand{\hga}{\hat{\cal G}({\cal A})}
\newcommand{\Edo}{\mbox{End}\:}
\newcommand{\for}{\mbox{for}}
\newcommand{\kn}{\mbox{ker}}
\newcommand{\Dlt}{\Delta}
\newcommand{\rad}{\mbox{Rad}}
\newcommand{\rta}{\rightarrow}
\newcommand{\mbb}{\mathbb}
\newcommand{\lra}{\Longrightarrow}
\newcommand{\X}{{\cal X}}
\newcommand{\Y}{{\cal Y}}
\newcommand{\Z}{{\cal Z}}
\newcommand{\U}{{\cal U}}
\newcommand{\V}{{\cal V}}
\newcommand{\W}{{\cal W}}
\setlength{\unitlength}{3pt}

\begin{center}{\Large \bf Polynomial Representation of $E_7$ and Its}\end{center}
\begin{center}{\Large \bf  Combinatorial and PDE Implications}
\end{center}
\vspace{0.2cm}

\begin{center}{\large Xiaoping Xu}\end{center}
\begin{center}{Institute of Mathematics, Academy of Mathematics \& System Sciences}\end{center}
\begin{center}{Chinese Academy of Sciences, Beijing 100190, P.R. China}
\footnote{Research supported
 by China NSF 10871193}\end{center}

\vspace{0.6cm}

 \begin{center}{\Large\bf Abstract}\end{center}

\vspace{1cm} {\small  In this paper, we use partial differential
equations to find the decomposition of the polynomial algebra over
the basic irreducible module of $E_7$ into a sum of irreducible
submodules.  Moreover, we obtain a combinatorial identity, saying
that the dimensions of certain irreducible modules of $E_7$ are
correlated by the binomial coefficients of fifty-five. Furthermore,
we prove that two families of irreducible submodules with three
integral parameters  are solutions of the fundamental invariant
differential operator corresponding to  Cartan's unique quartic
$E_7$ invariant.}

\section{Introduction}

The $E_7$ Lie algebra and group are important mathematical objects
with broad applications. They are second most complicated
finite-dimensional simple Lie algebra and group. The complexity
indeed imply rich connotation. Here we are only be able to list a
small part of it. Ramond [32] (1977) gave a group theoretical
analysis of a symmetry breaking affected by Higgs fields for the
vector-like unified theory base on $E_7$. Cvitanovi\'{e} [9]
(1981) studied the $E_7$ symmetry with negative bosonic dimension.
Han, Kim and Tanii [21] (1986) presented a supersymmetrization of
the six-dimensional anomaly-free $E_6\times E_7\times U(1)$ theory
with Lorentz Chern-Simons term. Kato and Kitazawa [23] (1989)
studied the correlation functions in an $E_7$-type
modular-invariant Wess-Zumino-Witten theory, which is related to
the scheme of string compatification proposed by Gepner. Ma [29]
(1990) found the spectrum-dependent solutions to the Yang-Baxter
equation for quantum $E_7$.

Fern\'{a}dez,  Garcia Fuertes and  Perelomov [16] (2005)
re-expressed the quantum Calogero-Sutherland model for the Lie
algebra $E_7$ and the particular value of the coupling constant
$k=1$, using the fundamental irreducible characters of the algebra
as dynamic variables. Moreover, they used the model to obtain
explicitly the characters and Clebsch-Gordan series for the algebra
$E_7$. D'Auria, Ferrara and Trigiante [10] (2006) considered
M-theory compactified on seven-torus with fluxes when all the seven
anti-symmetric tensor fields in four dimensions have been dualized
into scalars and thus the $E_{7(7)}$ symmetry was recovered. Duff
and  Ferrara [14] (2007) proposed that a particular tripartite
entanglement of seven quits, encoded in the Fano plane, is described
by the $E_7$ group and that the entanglement measure is given by
Cartan's quartic $E_7$ invariant. Brink,  Kim and  Ramond [1] (2008)
used the Cremmer-Julia $E_{7(7)}$ nonlinear symmetry of $N=8$
Supergravity to derive its order $\kappa$ on-shell Hamltonian in
terms of the chiral light-cone superfield. Borsten [5] (2008) used
the Freudenthal triple system to sketch  the precise dictionary
relating the 56 charges, parametrizing the general black hole
solution, to the 56-dimensional quantum states of $E_7$.

Part of the mathematical story of $E_7$ is as follows. Brown [2]
(1968) constructed a central simple nonassociative algebraic
structure on the Freudenthal triple system of $E_7$ and obtained a
new interpretation of the $E_7$ Lie algebra. He [3] (1969) proved
that the group leaving Cartan's quartic form invariant modulo its
center is exactly a simple $E_7$ Lie group. Faulkner [15] (1972)
studied the geometry of the planes of the octave symplectic
geometry of the ternary algebra related to the minimal
representation of $E_7$. Ferrar [17] (1980) used the Freudenthal
triple system  to classify the simple $E_7$ Lie algebras over
algebraic number fields. Kleidman and Ryba [24] (1993) proved a
Kostant's conjecture in the case of $E_7$, which says that the
group $PSL(2,37)$ embeds in an $E_7$ group. Moreover, Griess and
Ryba [20] (1994) proved that the group $U_3(8)$ embeds in an $E_7$
group and the group $Sz(8)$ does not. Ginzburg [19] (1995) showed
that the twisted partial $L$-function on the 56-dimensional
representation of $GE_7(\mbb{C})$ is entire except the points 0
and 1.

Cooperstein [8] (1995) proved that  Cartan's  quartic
$E_7$-invariant is the unique fundamental invariant over the basic
module. Moreover, Shult [34] (1997) showed that the basic
representation of $E_7$ provides the absolutely universal embedding
of the point-line geometry $E_{7,1}$. Plotkin [31] (1998)
investigated the stability of the $K_1$-functor for the $E_7$ group.
Dokovi\'{c} [12] (1999) classified nilpotent adjoint orbits of real
simple noncompact groups of type $E_7$ by means of Caley triples. He
[13] (2001) also worked out the partial order on nilpotent orbits in
the split real Lie algebra of type $E_7$. The $E_7$ root system was
used by Sekiguchi [33] (1999) to study the configurations of seven
lines on the real projective plane. Kleidman, Meierfrankenfeld and
Ryba [25] (1999) proved that the group $HS$ embeds in $E_7(5)$. Choi
and Yoon [6] (1999) calculated the homology of the double and triple
loop spaces of $E_7$.

Ukai [35] (2001) computed the $b$-functions of the prehomogeneous
vector space arising from a cuspidal character sheaf of $E_7$.
Garibaidi [18] (2001) gave  an explicit description of the
homogeneous projective varieties associated to the $E_7$ group with
trivial Tits algebras. A characterization of quadratic forms of type
$E_7$ was given by De Medits [11] (2002). Cherenousov [7] (2003)
used Hasse principle to prove that the Rost invariant has trivial
kernel for a quasi-split $E_7$ group. Kono, Lin and Nishimura [26]
(2003) characterized the mod 3 cohomology of $E_7$. Weiss [36]
(2006) solved a fundamental question about the structure of the
automorphism group of the Moufang quadrangles of type $E_7$.

In [43] and [44], we used partial differential equations to
decompose the polynomial algebra over the basic irreducible modules
of $F_4$ and $E_6$ into a sum of irreducible submodules.
Consequently, we obtained combinatorial identities of correlating
certain dimensions of irreducible modules. The motivation of this
work is to understand the functional impact of the $E_7$ Lie algebra
(group). We  use partial differential equations and the result in
[44] to decompose the polynomial algebra over the basic irreducible
module of $E_7$ into a sum of irreducible submodules. As a result,
we obtain a combinatorial identity, which  says that the dimensions
of certain irreducible modules of $E_7$ are correlated by the
binomial coefficients of fifty-five. Furthermore, we find that two
families of irreducible submodules with three integral parameters
are solutions of the fundamental invariant differential operator
corresponding to Cartan's unique quartic $E_7$-invariant. Below we
give a more detailed introduction to our results.

It has been known for may years that the representation theory of
Lie algebras is closely related to combinatorial identities.
Macdonald [30] generalized the Weyl denominator identities for
finite root systems to those for infinite affine root systems, which
are now known as the Macdonald's identities. We present here a
consequence  of the Macdonald's identities taken from Kostant's work
[27]. Let ${\cal G}$ be a finite-dimensional simple Lie algebra over
the field $\mbb{C}$ of complex numbers. Denote by $\Lmd^+$ the set
of dominant integral weights of ${\cal G}$ and by $V(\lmd)$ the
finite-dimensional irreducible ${\cal G}$-module with highest weight
$\lmd$. It is known that the Casimir operator takes a constant
$c(\lmd)$ on $V(\lmd)$. Macdonald's Theorem implies that there
exists a map $\chi:\Lmd^+\rta\{-1,0,1\}$ such that the following
identity holds:
$$(\prod_{n=1}^\infty(1-q^n))^{\dim {\cal
G}}=\sum_{\lmd\in\Lmd^+}\chi(\lmd)(\dim
V(\lmd))q^{c(\lmd)}.\eqno(1.1)$$ Kostant [28]  connected the above
identity to the abelian subalgebras of ${\cal G}$.

Our idea of using partial differential equations to solve Lie
algebra problems started in our earlier works [37] and [38] when we
tried to find functional generators for the invariants over
curvature tensor fields and for the differential invariants of
classical groups. Later we used partial differential equations to
find the explicit formulas for all the singular vectors in the Verma
modules of $sl(n,\mbb{C})$ (cf. [39]) and $sp(4,\mbb{C})$ (cf.
[40]). We [41] developed the methods of solving linear partial
differential equations of flag type. These equations naturally
appear in the problem of determining irreducible highest-weight
submodules of a general module for a finite-dimensional simple Lie
algebra, such as polynomial algebras over irreducible modules. Many
important linear partial differential equations in physics and
geometry are of flag type. In particular, we find new special
functions by which we are able to explicitly give the solutions of
the initial value problems of a large family of constant-coefficient
linear partial differential equations in terms of their
coefficients.

The Dynkin diagram of $E_7$ is as follows:

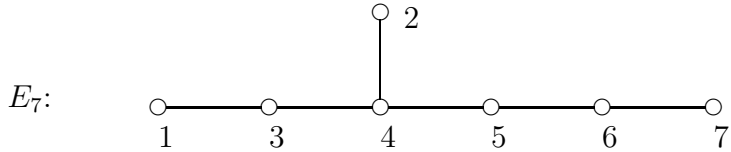
\begin{picture}(93,20)
\put(2,0){$E_7$:}\put(21,0){\circle{2}}\put(21,
-5){1}\put(22,0){\line(1,0){12}}\put(35,0){\circle{2}}\put(35,
-5){3}\put(36,0){\line(1,0){12}}\put(49,0){\circle{2}}\put(49,
-5){4}\put(49,1){\line(0,1){10}}\put(49,12){\circle{2}}\put(52,10){2}\put(50,0){\line(1,0){12}}
\put(63,0){\circle{2}}\put(63,-5){5}\put(64,0){\line(1,0){12}}\put(77,0){\circle{2}}\put(77,
-5){6}\put(78,0){\line(1,0){12}}\put(91,0){\circle{2}}\put(91,
-5){7}
\end{picture}
\vspace{0.7cm}

\noindent Denote by $\lmd_i$ the $i$th fundamental weight of $E_7$
with respect to the above labeling. Let $V$ be the 56-dimensional
irreducible $E_7$-module of highest weight $\lmd_7$. Denote by
${\cal A}$ the polynomial algebra (equivalently, symmetric tensor)
over $V$. A {\it singular vector} in ${\cal A}$ is a weight vector
annihilated by positive root vectors. The following is the main
theorem of this paper:\psp

{\bf Main Theorem}: {\it Any singular vector is a polynomial in a
linear singular vector $x_1$ of weight $\lmd_7$, a quadratic
singular vector $\zeta_1$ of weight $\lmd_1$, a cubic singular
vector $\vt$ of weight $\lmd_7$, a quartic singular vector $\vs$ of
weight $\lmd_6$ and a quartic singular vector $\eta$ of weight 0
(Cartan's invariant). Let $L(n_1,n_2,n_3,n_4,n_5)$ be the
irreducible submodule generated by
$\zeta_1^{n_1}\vs^{n_2}x_1^{n_3}\vt^{n_4}\eta^{n_5}$ with highest
weight $n_1\lmd_1+n_2\lmd_6+(n_3+n_4)\lmd_7$. Then
$${\cal A}=\bigoplus_{n_1,n_2,n_3,n_4,n_5=0}^\infty
L(n_1,n_2,n_3,n_4,n_5).\eqno(1.2)$$ In particular,
\begin{eqnarray*}& &(1-q)^{55}\sum_{n_1,n_2,n_3,n_4=0}^\infty(\mbox{dim}\:
V(n_1\lmd_1+n_2\lmd_6+(n_3+n_4)\lmd_7))q^{2n_1+4n_2+n_3+3n_4}\\
&=&1+q+q^2+q^3.\hspace{11.3cm}(1.3)\end{eqnarray*}
  Let ${\cal D}$ be the invariant differential operator dual to $\eta$. Then
  $$\sum_{n_1,n_2,n_3=0}^\infty(L(n_1,n_2,n_3,0,0)+L(n_1,n_2,n_3,1,0))\subset
  \{f\in{\cal A}\mid {\cal D}(f)=0\}.\eqno(1.4)$$}

\pse

Note that the identities (1.1) and (1.3) are dimensional properties
of irreducible modules. Our identity (1.3) says that the numbers
\begin{eqnarray*}& &\dim
V(n_1\lmd_1+n_2\lmd_6+n_3\lmd_7)\\&=&\frac{(n_2+4)(n_2+n_3+5)[\prod_{i=1}^7(n_1+i)(n_2+i)][\prod_{r=2}^8(n_2+n_3+r)]
[\prod_{s=5}^{11}(n_1+n_2+s)]}{2^3.3^3.4^4.5^5.6^5.7^5.8^4.9^4.10^3.11^3.12^2.13^2.14.15.16.17}
\\ & &\times[\prod_{p=6}^{12}(n_1+n_2+n_3+p)]
[\prod_{\iota=10}^{16}(n_1+2n_2+n_3+\iota)]
(2n_2+n_3+9)(n_1+n_2+8)\\&
&\times(n_1+n_2+n_3+9)(n_1+2n_2+n_3+13)(2n_1+2n_2+n_3+17)(n_3+1)\hspace{1.2cm}(1.17)
\end{eqnarray*}
 are correlated by the
binomial coefficients of fifty-five.

In Section 2, we explicitly construct the 56-dimensional basic
representation of $E_7$ in terms of differential operators  via the
root lattice construction of the $E_8$ simple Lie algebra. Explicit
constructions of the singular vectors $\zeta_1,\vs,\vt$ and $\eta$
are given in Section 3. The proof of the main theorem is presented
in Section 4.

 \section{Basic Representation of $E_7$}

In this section, we will explicitly construct the 56-dimensional
basic irreducible representation of $E_7$.

 For convenience, we will use the notion
$$\ol{i,i+j}=\{i,i+1,i+2,...,i+j\}\eqno(2.1)$$
for integer $i$ and positive integer $j$ throughout this paper. We
start with the root lattice construction of the simple Lie algebra
of type $E_8$. As we all known, the Dynkin diagram of $E_8$ is as
follows:

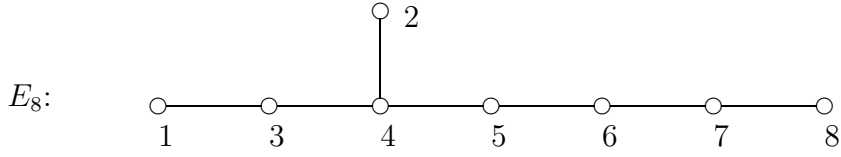
\begin{picture}(110,23)
\put(2,0){$E_8$:}\put(21,0){\circle{2}}\put(21,
-5){1}\put(22,0){\line(1,0){12}}\put(35,0){\circle{2}}\put(35,
-5){3}\put(36,0){\line(1,0){12}}\put(49,0){\circle{2}}\put(49,
-5){4}\put(49,1){\line(0,1){10}}\put(49,12){\circle{2}}\put(52,10){2}\put(50,0){\line(1,0){12}}
\put(63,0){\circle{2}}\put(63,-5){5}\put(64,0){\line(1,0){12}}\put(77,0){\circle{2}}\put(77,
-5){6}\put(78,0){\line(1,0){12}}\put(91,0){\circle{2}}\put(91,
-5){7}\put(92,0){\line(1,0){12}}\put(105,0){\circle{2}}\put(105,
-5){8}
\end{picture}
\vspace{0.7cm}

 \noindent Let $\{\al_i\mid i\in\ol{1,8}\}$ be the
simple positive roots corresponding to the vertices in the diagram,
and let $\Phi_{E_8}$ be the root system of $E_8$. Set
$$Q_{E_8}=\sum_{i=1}^8\mbb{Z}\al_i,\eqno(2.2)$$ the root lattice of type
$E_8$. Denote by $(\cdot,\cdot)$ the symmetric $\mbb{Z}$-bilinear
form on $Q_{E_8}$ such that
$$\Phi_{E_8}=\{\al\in Q_{E_8}\mid (\al,\al)=2\}.\eqno(2.3)$$
Define $F(\cdot,\cdot):\; Q_{E_8}\times  Q_{E_8}\rta \{\pm 1\}$ by
$$F(\sum_{i=1}^8k_i\al_i,\sum_{j=1}^8l_j\al_j)=(-1)^{\sum_{i=1}^8k_il_i+\sum_{8\geq i>j\geq 1}k_il_j
(\al_i,\al_j)},\qquad k_i,l_j\in\mbb{Z}.\eqno(2.4)$$ Then for
$\al,\be,\gm\in  Q_{E_8}$,
$$F(\al+\be,\gm)=F(\al,\gm)F(\be,\gm),\;\;F(\al,\be+\gm)=F(\al,\be)F(\al,\gm),\eqno(2.5)
$$
$$F(\al,\be)F(\be,\al)^{-1}=(-1)^{(\al,\be)},\;\;F(\al,\al)=(-1)^{(\al,\al)/2}.
\eqno(2.6)$$ In particular,
$$F(\al,\be)=-F(\be,\al)\qquad
\mbox{if}\;\;\al,\be,\al+\be\in \Phi_{E_8}.\eqno(2.7)$$

Denote
$$H_{E_8}=\sum_{i=1}^8\mbb{C}\al_i.\eqno(2.8)$$
The simple Lie algebra of type $E_8$ is
$${\cal
G}^{E_8}=H_{E_8}\oplus\bigoplus_{\al\in
\Phi_{E_8}}\mbb{C}E_{\al}\eqno(2.9)$$ with the Lie bracket
$[\cdot,\cdot]$ determined by:
 $$[H_{E_8},H_{E_8}]=0,\;\;[h,E_{\al}]=(h,\al)E_{\al},\;\;[E_{\al},E_{-\al}]=-\al,
 \eqno(2.10)$$
 $$[E_{\al},E_{\be}]=\left\{\begin{array}{ll}0&\mbox{if}\;\al+\be\not\in \Phi_{E_8},\\
 F(\al,\be)E_{\al+\be}&\mbox{if}\;\al+\be\in\Phi_{E_8}.\end{array}\right.\eqno(2.11)$$
for $\al,\be\in\Phi_{E_8}$ and $h\in H_{E_8}$.

Note that
$$Q_{E_7}=\sum_{i=1}^7\mbb{Z}\al_i\subset Q_{E_8}\eqno(2.12)$$
is the root lattice of $E_7$ and
$$\Phi_{E_7}=Q_{E_7}\bigcap \Phi_{E_8}\eqno(2.13)$$
is the root system of $E_7$. Set
$$H_{E_7}=\sum_{i=1}^7\mbb{C}\al_i.\eqno(2.14)$$
Then the subalgebra
$${\cal G}^{E_7}=H_{E_7}\oplus\bigoplus_{\al\in
\Phi_{E_7}}\mbb{C}E_{\al}\eqno(2.15)$$ of ${\cal G}^{E_8}$ is
exactly the simple Lie algebra of ${\cal G}^{E_7}$. Denote by
$\Phi_{E_7}^+$ the set of positive roots of $E_7$ and by
$\Phi_{E_8}^+$ the set of positive roots of $E_8$. The elements of
$\Phi_{E_7}^+$ are:
$$\al_1+\al_2+2\al_3+3\al_4+2\al_5+\al_6+\sum_{r=1}^7\al_r,\eqno(2.15)$$
$$\{\al_1+2\al_2+2\al_3+3\al_4+2\al_5+\al_6+\sum_{r=i+1}^7\al_r
\mid i\in\ol{2,7}\},\eqno(2.16)$$
$$\{\al_1+\sum_{r=3}^j\al_r\mid j\in\ol{2,7}\}\bigcup
\{\sum_{r=i+1}^j\al_r\mid 2\leq i<j\leq 7\},\eqno(2.17)$$
$$\{\sum_{s=2}^j\al_s+\sum_{t=4}^k\al_t\mid 2\leq j<k\leq
7\},\eqno(2.18)$$
$$\{\sum_{\iota=1}^i\al_\iota+
\sum_{s=3}^j\al_s+\sum_{t=4}^k\al_t\mid 2\leq i< j<k\leq
7\}.\eqno(2.19)$$

Denote by $\bar\Phi_{E_8}^+$ the set of the following positive
roots:
$$\al_1+\sum_{r=3}^8\al_r,\qquad\sum_{r=i+1}^8\al_r\qquad\for\;\;i\in\ol{2,7},\eqno(2.20)$$
$$\al_4+\sum_{s=3}^5\al_s+\sum_{p=1}^6\al_p+\sum_{r=1}^7\al_r+\sum_{s=j+1}^8\al_s
\qquad\for\;\;j\in\ol{2,7},\eqno(2.21)$$
$$\al_4+\sum_{s=2}^5\al_s+\sum_{p=1}^6\al_p+\sum_{r=i+1}^7\al_r+\sum_{s=j+1}^8\al_s
\qquad\for\;\;2\leq i<j\leq 7,\eqno(2.22)$$
$$\sum_{s=2}^j\al_s+\sum_{t=4}^8\al_t\qquad\for\;\;j\in\ol{2,7},\eqno(2.23)$$
$$\sum_{\iota=1}^i\al_\iota+
\sum_{s=3}^j\al_s+\sum_{t=4}^8\al_t\qquad\for\;\;2\leq i< j\leq 7,
\eqno(2.24)$$
$$\al_1+3\al_2+3\al_3+5\al_4+4\al_5+3\al_6+2\al_7+\al_8,\eqno(2.25)$$
$$2\al_1+2\al_2+3\al_3+5\al_4+4\al_5+3\al_6+2\al_7+\al_8+
\sum_{r=2}^i\al_r\;\;\for\;\;i\in\ol{2,8}.\eqno(2.26)$$ Then
$$\Phi_{E_8}^+=\Phi_{E_7}^+\bigcup \bar\Phi_{E_8}^+.\eqno(2.27)$$
Set
$$\theta=2\al_1+2\al_2+3\al_3+5\al_4+4\al_5+3\al_6+2\al_7+\al_8+
\sum_{r=2}^8\al_r.\eqno(2.28)$$ In particular,
$$V=\sum_{\theta\neq\be\in \bar\Phi_{E_8}^+}\mbb{C}E_\be\eqno(2.28)$$
forms the 56-dimensional basic ${\cal G}^{E_7}$-module of highest
weight $\lmd_7$ with the representation $\mbox{ad}_{{\cal
G}^{E_8}}$.

For convenience, we also denote
$$E_{(k_1,...,k_r)}=E_\al\qquad\for\;\;\al=\sum_{s=1}^rk_s\al_s\in\Phi_{E_8}^+,\;k_r\neq
0.\eqno(2.29)$$ Let
$$x_1=E_{(2,3,4,6,5,4,3,1)},\qquad x_2=E_{(2,3,4,6,5,4,2,1)},\qquad
x_3=E_{(2,3,4,6,5,3,2,1)},\eqno(2.30)$$
$$x_4=E_{(2,3,4,6,4,3,2,1)},\qquad x_5=E_{(2,3,4,5,4,3,2,1)},\qquad
x_6=E_{(2,3,3,5,4,3,2,1)},\eqno(2.31)$$
$$x_7=E_{(2,2,4,5,4,3,2,1)},\qquad x_8=E_{(1,3,3,5,4,3,2,1)},\qquad
x_9=E_{(2,2,3,5,4,3,2,1)},\eqno(2.32)$$
$$x_{10}=E_{(2,2,3,4,4,3,2,1)},\qquad x_{11}=E_{(1,2,3,5,4,3,2,1)},\qquad
x_{12}=E_{(2,2,3,4,3,3,2,1)}, \eqno(2.33)$$
$$ x_{13}=E_{(1,2,3,4,4,3,2,1)},\qquad x_{14}=E_{(2,2,3,4,3,2,2,1)},\qquad
x_{15}=E_{(1,2,2,4,4,3,2,1)},\eqno(2.34)$$
$$x_{16}=E_{(1,2,3,4,3,3,2,1)},\qquad
x_{17}=E_{(2,2,3,4,3,2,1,1)},\qquad
x_{18}=E_{(1,2,2,4,3,3,2,1)},\eqno(2.35)$$
$$x_{19}=E_{(1,2,3,4,3,2,2,1)},\qquad
x_{20}=E_{(1,2,2,3,3,3,2,1)},\qquad
x_{21}=E_{(1,2,2,4,3,2,2,1)},\eqno(2.36)$$
$$x_{22}=E_{(1,2,3,4,3,2,1,1)},
\qquad x_{23}=E_{(1,1,2,3,3,3,2,1)},\qquad
x_{24}=E_{(1,2,2,3,3,2,2,1)}, \eqno(2.37)$$
$$x_{25}=E_{(1,2,2,4,3,2,1,1)},\qquad
x_{26}=E_{(1,1,2,3,3,2,2,1)},\qquad x_{27}=E_{(1,2,2,3,2,2,2,1)},
\eqno(2.38)$$
$$x_{28}=E_{(1,2,2,3,3,2,1,1)},\qquad
x_{29}=E_{(1,1,2,3,2,2,2,1)},\qquad x_{30}=E_{(1,1,2,3,3,2,1,1)},
 \eqno(2.39)$$
$$x_{31}=E_{(1,2,2,3,2,2,1,1)},\qquad x_{32}=E_{(1,1,2,2,2,2,2,1)},
\qquad x_{33}=E_{(1,1,2,3,2,2,1,1)}, \eqno(2.40)$$
$$ x_{34}=E_{(1,2,2,3,2,1,1,1)},\qquad
x_{35}=E_{(1,1,1,2,2,2,2,1)},\qquad
x_{36}=E_{(1,1,2,2,2,2,1,1)},\eqno(2.41)$$
$$x_{37}=E_{(1,1,2,3,2,1,1,1)},\qquad x_{38}=E_{(1,1,1,2,2,2,1,1)},\qquad x_{39}=E_{(1,1,2,2,2,1,1,1)},
\eqno(2.45)$$
$$x_{40}=E_{(0,1,1,2,2,2,2,1)}, \qquad x_{41}=E_{(1,1,1,2,2,1,1,1)},
\qquad x_{42}=E_{(1,1,2,2,1,1,1,1)},\eqno(2.46)$$
$$x_{43}=E_{(0,1,1,2,2,2,1,1)},\qquad
x_{44}=E_{(1,1,1,2,1,1,1,1)},\qquad
x_{45}=E_{(0,1,1,2,2,1,1,1)},\eqno(2.47)$$
$$x_{46}=E_{(1,1,1,1,1,1,1,1)},\qquad x_{47}=E_{(0,1,1,2,1,1,1,1)},\qquad
x_{48}=E_{(0,1,1,1,1,1,1,1)},\eqno(2.48)$$
$$x_{49}=E_{(1,0,1,1,1,1,1,1)},\qquad
x_{50}=E_{(0,1,0,1,1,1,1,1)},\qquad
x_{51}=E_{(0,0,1,1,1,1,1,1)},\eqno(2.49)$$
$$x_{52}=E_{(0,0,0,1,1,1,1,1)},\qquad
x_{53}=E_{(0,0,0,0,1,1,1,1)},\qquad
x_{54}=E_{(0,0,0,0,0,1,1,1)},\eqno(2.50)$$
$$x_{55}=E_{(0,0,0,0,0,0,1,1)},\qquad
x_{56}=E_{(0,0,0,0,0,0,0,1)}.\eqno(2.51)$$ Then $\{x_i\mid
i\in\ol{1,56}\}$ forms a basis of $V$. With respect to the basis, we
have the following representation formulas:
\begin{eqnarray*}E_{\al_1}|_V&=&-x_6\ptl_{x_8}-x_9\ptl_{x_{11}}
-x_{10}\ptl_{x_{13}}-x_{12}\ptl_{x_{16}}-x_{14}\ptl_{x_{19}}-x_{17}\ptl_{x_{22}}
\\ &&+x_{35}\ptl_{x_{40}}+x_{38}\ptl_{x_{43}}+x_{41}\ptl_{x_{45}}
+x_{44}\ptl_{x_{47}}+x_{46}\ptl_{x_{48}}+x_{49}\ptl_{x_{51}},\hspace{2.4cm}(2.52)
\end{eqnarray*}
\begin{eqnarray*}E_{\al_2}|_V&=&x_5\ptl_{x_7}+x_6\ptl_{x_9}
+x_8\ptl_{x_{11}}-x_{20}\ptl_{x_{23}}-x_{24}\ptl_{x_{26}}-x_{27}\ptl_{x_{29}}
\\ &&-x_{28}\ptl_{x_{30}}-x_{31}\ptl_{x_{33}}-x_{34}\ptl_{x_{37}}
+x_{46}\ptl_{x_{49}}+x_{48}\ptl_{x_{51}}+x_{50}\ptl_{x_{52}},\hspace{2.4cm}(2.53)
\end{eqnarray*}
\begin{eqnarray*}E_{\al_3}|_V&=&-x_5\ptl_{x_6}-x_7\ptl_{x_9}
-x_{13}\ptl_{x_{15}}-x_{16}\ptl_{x_{18}}-x_{19}\ptl_{x_{21}}-x_{22}\ptl_{x_{25}}
\\ &&+x_{32}\ptl_{x_{35}}+x_{36}\ptl_{x_{38}}+x_{39}\ptl_{x_{41}}
+x_{42}\ptl_{x_{44}}+x_{48}\ptl_{x_{50}}+x_{51}\ptl_{x_{52}},
\hspace{2.4cm}(2.54)
\end{eqnarray*}
\begin{eqnarray*}E_{\al_4}|_V&=&x_4\ptl_{x_5}-x_9\ptl_{x_{10}}
-x_{11}\ptl_{x_{13}}-x_{18}\ptl_{x_{20}}-x_{21}\ptl_{x_{24}}-x_{25}\ptl_{x_{28}}
\\ &&-x_{29}\ptl_{x_{32}}-x_{33}\ptl_{x_{36}}-x_{37}\ptl_{x_{39}}
-x_{44}\ptl_{x_{46}}-x_{47}\ptl_{x_{48}}+x_{52}\ptl_{x_{53}},
\hspace{2.4cm}(2.55)
\end{eqnarray*}
\begin{eqnarray*}E_{\al_5}|_V&=&x_3\ptl_{x_4}-x_{10}\ptl_{x_{12}}
-x_{13}\ptl_{x_{16}}-x_{15}\ptl_{x_{18}}-x_{24}\ptl_{x_{27}}-x_{26}\ptl_{x_{29}}
\\ &&-x_{28}\ptl_{x_{31}}-x_{30}\ptl_{x_{33}}-x_{39}\ptl_{x_{42}}
-x_{41}\ptl_{x_{44}}-x_{45}\ptl_{x_{47}}+x_{53}\ptl_{x_{54}},
\hspace{2.4cm}(2.56)
\end{eqnarray*}
\begin{eqnarray*}E_{\al_6}|_V&=&x_2\ptl_{x_3}-x_{12}\ptl_{x_{14}}
-x_{16}\ptl_{x_{19}}-x_{18}\ptl_{x_{21}}-x_{20}\ptl_{x_{24}}-x_{23}\ptl_{x_{26}}
\\ &&-x_{31}\ptl_{x_{34}}-x_{33}\ptl_{x_{37}}-x_{36}\ptl_{x_{39}}
-x_{38}\ptl_{x_{41}}-x_{43}\ptl_{x_{45}}+x_{54}\ptl_{x_{55}},
\hspace{2.4cm}(2.57)
\end{eqnarray*}
\begin{eqnarray*}E_{\al_7}|_V&=&x_1\ptl_{x_2}-x_{14}\ptl_{x_{17}}
-x_{19}\ptl_{x_{22}}-x_{21}\ptl_{x_{25}}-x_{24}\ptl_{x_{28}}-x_{26}\ptl_{x_{30}}
\\ &&-x_{27}\ptl_{x_{31}}-x_{29}\ptl_{x_{33}}-x_{32}\ptl_{x_{36}}
-x_{35}\ptl_{x_{38}}-x_{40}\ptl_{x_{43}}+x_{55}\ptl_{x_{56}},
\hspace{2.4cm}(2.58)
\end{eqnarray*}
\begin{eqnarray*}E_{(1,0,1)}|_V&=&-x_5\ptl_{x_8}-x_7\ptl_{x_{11}}
+x_{10}\ptl_{x_{15}}+x_{12}\ptl_{x_{18}}+x_{14}\ptl_{x_{21}}+x_{17}\ptl_{x_{25}}
\\ &&-x_{32}\ptl_{x_{40}}-x_{36}\ptl_{x_{43}}-x_{39}\ptl_{x_{45}}
-x_{42}\ptl_{x_{47}}+x_{46}\ptl_{x_{50}}+x_{49}\ptl_{x_{52}},
\hspace{1.9cm}(2.59)
\end{eqnarray*}
\begin{eqnarray*}E_{(0,1,0,1)}|_V&=&-x_4\ptl_{x_7}-x_6\ptl_{x_{10}}
-x_8\ptl_{x_{13}}-x_{18}\ptl_{x_{23}}-x_{21}\ptl_{x_{26}}-x_{25}\ptl_{x_{30}}
\\ &&+x_{27}\ptl_{x_{32}}+x_{31}\ptl_{x_{36}}+x_{34}\ptl_{x_{39}}
+x_{44}\ptl_{x_{49}}+x_{47}\ptl_{x_{51}}+x_{50}\ptl_{x_{53}},
\hspace{1.6cm}(2.60)
\end{eqnarray*}
\begin{eqnarray*}E_{(0,0,1,1)}|_V&=&x_4\ptl_{x_6}+x_7\ptl_{x_{10}}
-x_{11}\ptl_{x_{15}}+x_{16}\ptl_{x_{20}}+x_{19}\ptl_{x_{24}}+x_{22}\ptl_{x_{28}}
\\ &&+x_{29}\ptl_{x_{35}}+x_{33}\ptl_{x_{38}}+x_{37}\ptl_{x_{41}}
-x_{42}\ptl_{x_{46}}+x_{47}\ptl_{x_{50}}+x_{51}\ptl_{x_{53}},
\hspace{1.6cm}(2.61)
\end{eqnarray*}
\begin{eqnarray*}& &E_{(0,0,0,1,1)}|_V=-x_3\ptl_{x_5}+x_9\ptl_{x_{12}}
+x_{11}\ptl_{x_{16}}-x_{15}\ptl_{x_{20}}+x_{21}\ptl_{x_{27}}+x_{25}\ptl_{x_{31}}
\\ &&-x_{26}\ptl_{x_{32}}-x_{30}\ptl_{x_{36}}+x_{37}\ptl_{x_{42}}
-x_{41}\ptl_{x_{46}}-x_{45}\ptl_{x_{48}}+x_{52}\ptl_{x_{54}},
\hspace{3.7cm}(2.62)
\end{eqnarray*}
\begin{eqnarray*}& &E_{(0,0,0,0,1,1)}|_V=-x_2\ptl_{x_4}+x_{10}\ptl_{x_{14}}
+x_{13}\ptl_{x_{19}}+x_{15}\ptl_{x_{21}}-x_{20}\ptl_{x_{27}}-x_{23}\ptl_{x_{29}}
\\ &&+x_{28}\ptl_{x_{34}}+x_{30}\ptl_{x_{37}}-x_{36}\ptl_{x_{42}}
-x_{38}\ptl_{x_{44}}-x_{43}\ptl_{x_{47}}+x_{53}\ptl_{x_{55}},
\hspace{3.7cm}(2.63)
\end{eqnarray*}
\begin{eqnarray*}& &E_{(0,0,0,0,0,1,1)}|_V=-x_1\ptl_{x_3}+x_{12}\ptl_{x_{17}}
+x_{16}\ptl_{x_{22}}+x_{18}\ptl_{x_{25}}+x_{20}\ptl_{x_{28}}+x_{23}\ptl_{x_{30}}
\\ &&-x_{27}\ptl_{x_{34}}-x_{29}\ptl_{x_{37}}-x_{32}\ptl_{x_{39}}
-x_{35}\ptl_{x_{41}}-x_{40}\ptl_{x_{45}}+x_{54}\ptl_{x_{56}},
\hspace{3.7cm}(2.64)
\end{eqnarray*}
\begin{eqnarray*}E_{(1,0,1,1)}|_V&=&x_4\ptl_{x_8}+x_7\ptl_{x_{13}}
+x_9\ptl_{x_{15}}-x_{12}\ptl_{x_{20}}-x_{14}\ptl_{x_{24}}-x_{17}\ptl_{x_{28}}
\\ &&-x_{29}\ptl_{x_{40}}-x_{33}\ptl_{x_{43}}-x_{37}\ptl_{x_{45}}
+x_{42}\ptl_{x_{48}}+x_{44}\ptl_{x_{50}}+x_{49}\ptl_{x_{53}},
\hspace{1.6cm}(2.65)
\end{eqnarray*}
\begin{eqnarray*}E_{(0,1,1,1)}|_V&=&-x_4\ptl_{x_9}+x_5\ptl_{x_{10}}
-x_8\ptl_{x_{15}}+x_{16}\ptl_{x_{23}}+x_{19}\ptl_{x_{26}}+x_{22}\ptl_{x_{30}}
\\ &&-x_{27}\ptl_{x_{35}}-x_{31}\ptl_{x_{38}}-x_{34}\ptl_{x_{41}}
+x_{42}\ptl_{x_{49}}-x_{47}\ptl_{x_{52}}+x_{48}\ptl_{x_{53}},
\hspace{1.6cm}(2.66)
\end{eqnarray*}
\begin{eqnarray*}& &E_{(0,1,0,1,1)}|_V=x_3\ptl_{x_7}+x_6\ptl_{x_{12}}
+x_8\ptl_{x_{16}}-x_{15}\ptl_{x_{23}}+x_{21}\ptl_{x_{29}}+x_{24}\ptl_{x_{32}}
\\ &&+x_{25}\ptl_{x_{33}}+x_{28}\ptl_{x_{36}}-x_{34}\ptl_{x_{42}}
+x_{41}\ptl_{x_{49}}+x_{45}\ptl_{x_{51}}+x_{50}\ptl_{x_{54}},
\hspace{3.7cm}(2.67)
\end{eqnarray*}
\begin{eqnarray*}& &E_{(0,0,1,1,1)}|_V=-x_3\ptl_{x_6}-x_7\ptl_{x_{12}}
+x_{11}\ptl_{x_{18}}+x_{13}\ptl_{x_{20}}-x_{19}\ptl_{x_{27}}-x_{22}\ptl_{x_{31}}
\\ &&+x_{26}\ptl_{x_{35}}+x_{30}\ptl_{x_{38}}-x_{37}\ptl_{x_{44}}
-x_{39}\ptl_{x_{46}}+x_{45}\ptl_{x_{50}}+x_{51}\ptl_{x_{54}},
\hspace{3.7cm}(2.68)
\end{eqnarray*}
\begin{eqnarray*}& &E_{(0,0,0,1,1,1)}|_V=x_2\ptl_{x_5}-x_9\ptl_{x_{14}}
-x_{11}\ptl_{x_{19}}+x_{15}\ptl_{x_{24}}+x_{18}\ptl_{x_{28}}-x_{23}\ptl_{x_{32}}
\\ &&-x_{25}\ptl_{x_{34}}+x_{30}\ptl_{x_{39}}+x_{33}\ptl_{x_{42}}
-x_{38}\ptl_{x_{46}}-x_{43}\ptl_{x_{48}}+x_{52}\ptl_{x_{55}},
\hspace{3.7cm}(2.69)
\end{eqnarray*}
\begin{eqnarray*}& &E_{(0,0,0,0,1,1,1)}|_V=x_1\ptl_{x_4}-x_{10}\ptl_{x_{17}}
-x_{13}\ptl_{x_{22}}-x_{15}\ptl_{x_{25}}+x_{20}\ptl_{x_{31}}+x_{23}\ptl_{x_{33}}
\\ &&+x_{24}\ptl_{x_{34}}+x_{26}\ptl_{x_{37}}-x_{32}\ptl_{x_{42}}
-x_{35}\ptl_{x_{44}}-x_{40}\ptl_{x_{47}}+x_{53}\ptl_{x_{56}},
\hspace{3.7cm}(2.70)
\end{eqnarray*}
\begin{eqnarray*}E_{(1,1,1,1)}|_V&=&-x_4\ptl_{x_{11}}+x_5\ptl_{x_{13}}
+x_6\ptl_{x_{15}}-x_{12}\ptl_{x_{23}}-x_{14}\ptl_{x_{26}}-x_{17}\ptl_{x_{30}}
\\ &&+x_{27}\ptl_{x_{40}}+x_{31}\ptl_{x_{43}}+x_{34}\ptl_{x_{45}}
-x_{42}\ptl_{x_{51}}-x_{44}\ptl_{x_{52}}+x_{46}\ptl_{x_{53}},
\hspace{1.6cm}(2.71)
\end{eqnarray*}
\begin{eqnarray*}& &E_{(1,0,1,1,1)}|_V=-x_3\ptl_{x_8}-x_7\ptl_{x_{16}}
-x_9\ptl_{x_{18}}-x_{10}\ptl_{x_{20}}+x_{14}\ptl_{x_{27}}+x_{17}\ptl_{x_{31}}
\\ &&-x_{26}\ptl_{x_{40}}-x_{30}\ptl_{x_{43}}+x_{37}\ptl_{x_{47}}+x_{39}\ptl_{x_{48}}
+x_{41}\ptl_{x_{50}}+x_{49}\ptl_{x_{54}}, \hspace{3.7cm}(2.72)
\end{eqnarray*}
\begin{eqnarray*}& &E_{(0,1,1,1,1)}|_V=x_3\ptl_{x_9}-x_5\ptl_{x_{12}}
+x_8\ptl_{x_{18}}+x_{13}\ptl_{x_{23}}-x_{19}\ptl_{x_{29}}-x_{22}\ptl_{x_{33}}
\\ &&-x_{24}\ptl_{x_{35}}-x_{28}\ptl_{x_{38}}+x_{34}\ptl_{x_{44}}
+x_{39}\ptl_{x_{49}}-x_{45}\ptl_{x_{52}}+x_{48}\ptl_{x_{54}},
\hspace{3.7cm}(2.73)
\end{eqnarray*}
\begin{eqnarray*}& &E_{(0,1,0,1,1,1)}|_V=-x_2\ptl_{x_7}-x_6\ptl_{x_{14}}
-x_8\ptl_{x_{19}}+x_{15}\ptl_{x_{26}}+x_{18}\ptl_{x_{30}}+x_{20}\ptl_{x_{32}}
\\ &&-x_{25}\ptl_{x_{37}}-x_{28}\ptl_{x_{39}}-x_{31}\ptl_{x_{42}}
+x_{38}\ptl_{x_{49}}+x_{43}\ptl_{x_{51}}+x_{50}\ptl_{x_{55}},
\hspace{3.7cm}(2.74)
\end{eqnarray*}
\begin{eqnarray*}& &E_{(0,0,1,1,1,1)}|_V=x_2\ptl_{x_6}+x_7\ptl_{x_{14}}
-x_{11}\ptl_{x_{21}}-x_{13}\ptl_{x_{24}}-x_{16}\ptl_{x_{28}}+x_{23}\ptl_{x_{35}}
\\ &&+x_{22}\ptl_{x_{34}}-x_{30}\ptl_{x_{41}}-x_{33}\ptl_{x_{44}}
-x_{36}\ptl_{x_{46}}+x_{43}\ptl_{x_{50}}+x_{51}\ptl_{x_{55}},
\hspace{3.7cm}(2.75)
\end{eqnarray*}
\begin{eqnarray*}& &E_{(0,0,0,1,1,1,1)}|_V=-x_1\ptl_{x_5}+x_9\ptl_{x_{17}}
+x_{11}\ptl_{x_{22}}-x_{15}\ptl_{x_{28}}-x_{18}\ptl_{x_{31}}-x_{21}\ptl_{x_{34}}
\\ &&+x_{23}\ptl_{x_{36}}+x_{26}\ptl_{x_{39}}+x_{29}\ptl_{x_{42}}
-x_{35}\ptl_{x_{46}}-x_{40}\ptl_{x_{48}}+x_{52}\ptl_{x_{56}},
\hspace{3.7cm}(2.76)
\end{eqnarray*}
\begin{eqnarray*}& &E_{(1,1,1,1,1)}|_V=x_3\ptl_{x_{11}}-x_5\ptl_{x_{16}}
-x_6\ptl_{x_{18}}-x_{10}\ptl_{x_{23}}+x_{14}\ptl_{x_{29}}+x_{17}\ptl_{x_{33}}
\\ &&+x_{24}\ptl_{x_{40}}+x_{28}\ptl_{x_{43}}-x_{34}\ptl_{x_{47}}
-x_{39}\ptl_{x_{51}}-x_{41}\ptl_{x_{52}}+x_{46}\ptl_{x_{54}},
\hspace{3.7cm}(2.77)
\end{eqnarray*}
\begin{eqnarray*}& &E_{(1,0,1,1,1,1)}|_V=x_2\ptl_{x_8}+x_7\ptl_{x_{19}}
+x_9\ptl_{x_{21}}+x_{10}\ptl_{x_{24}}+x_{12}\ptl_{x_{28}}-x_{17}\ptl_{x_{34}}
\\ &&-x_{23}\ptl_{x_{40}}+x_{30}\ptl_{x_{45}}+x_{33}\ptl_{x_{47}}
+x_{36}\ptl_{x_{48}}+x_{38}\ptl_{x_{50}}+x_{49}\ptl_{x_{55}},
\hspace{3.7cm}(2.78)
\end{eqnarray*}
\begin{eqnarray*}& &E_{(0,1,1,2,1)}|_V=-x_3\ptl_{x_{10}}+x_4\ptl_{x_{12}}
-x_8\ptl_{x_{20}}+x_{11}\ptl_{x_{23}}+x_{19}\ptl_{x_{32}}-x_{21}\ptl_{x_{35}}
\\ &&+x_{22}\ptl_{x_{36}}-x_{25}\ptl_{x_{38}}-x_{34}\ptl_{x_{46}}
+x_{37}\ptl_{x_{49}}-x_{45}\ptl_{x_{53}}+x_{47}\ptl_{x_{54}},
\hspace{3.7cm}(2.79)
\end{eqnarray*}
\begin{eqnarray*}& &E_{(0,1,1,1,1,1)}|_V=-x_2\ptl_{x_9}+x_5\ptl_{x_{14}}
-x_8\ptl_{x_{21}}-x_{13}\ptl_{x_{26}}-x_{16}\ptl_{x_{30}}-x_{20}\ptl_{x_{35}}
\\ &&+x_{22}\ptl_{x_{37}}+x_{28}\ptl_{x_{41}}+x_{31}\ptl_{x_{44}}
+x_{36}\ptl_{x_{49}}-x_{43}\ptl_{x_{52}}+x_{48}\ptl_{x_{55}},
\hspace{3.7cm}(2.80)
\end{eqnarray*}
\begin{eqnarray*}& &E_{(0,1,0,1,1,1,1)}|_V=x_1\ptl_{x_7}+x_6\ptl_{x_{17}}
+x_8\ptl_{x_{22}}-x_{15}\ptl_{x_{30}}-x_{18}\ptl_{x_{33}}-x_{21}\ptl_{x_{37}}
\\ &&-x_{20}\ptl_{x_{36}}-x_{24}\ptl_{x_{39}}-x_{27}\ptl_{x_{42}}
+x_{35}\ptl_{x_{49}}+x_{40}\ptl_{x_{51}}+x_{50}\ptl_{x_{56}},
\hspace{3.7cm}(2.81)
\end{eqnarray*}
\begin{eqnarray*}& &E_{(0,0,1,1,1,1,1)}|_V=-x_1\ptl_{x_6}-x_7\ptl_{x_{17}}
+x_{11}\ptl_{x_{25}}+x_{13}\ptl_{x_{28}}+x_{16}\ptl_{x_{31}}+x_{19}\ptl_{x_{34}}
\\ &&-x_{23}\ptl_{x_{38}}-x_{26}\ptl_{x_{41}}-x_{29}\ptl_{x_{44}}
-x_{32}\ptl_{x_{46}}+x_{40}\ptl_{x_{50}}+x_{51}\ptl_{x_{56}},
\hspace{3.7cm}(2.82)
\end{eqnarray*}
\begin{eqnarray*}& &E_{(1,1,1,2,1)}|_V=-x_3\ptl_{x_{13}}+x_4\ptl_{x_{16}}
+x_6\ptl_{x_{20}}-x_9\ptl_{x_{23}}-x_{14}\ptl_{x_{32}}+x_{21}\ptl_{x_{40}}
\\ &&-x_{17}\ptl_{x_{36}}+x_{25}\ptl_{x_{43}}+x_{34}\ptl_{x_{48}}
-x_{37}\ptl_{x_{51}}-x_{41}\ptl_{x_{53}}+x_{44}\ptl_{x_{54}},
\hspace{3.7cm}(2.83)
\end{eqnarray*}
\begin{eqnarray*}& &E_{(1,1,1,1,1,1)}|_V=-x_2\ptl_{x_{11}}+x_5\ptl_{x_{19}}
+x_6\ptl_{x_{21}}+x_{10}\ptl_{x_{26}}+x_{12}\ptl_{x_{29}}-x_{17}\ptl_{x_{37}}
\\ &&+x_{20}\ptl_{x_{40}}-x_{28}\ptl_{x_{45}}-x_{31}\ptl_{x_{47}}
-x_{36}\ptl_{x_{51}}-x_{38}\ptl_{x_{52}}+x_{46}\ptl_{x_{55}},
\hspace{3.7cm}(2.84)
\end{eqnarray*}
\begin{eqnarray*}& &E_{(1,0,1,1,1,1,1)}|_V=-x_1\ptl_{x_8}-x_7\ptl_{x_{22}}
-x_9\ptl_{x_{25}}-x_{10}\ptl_{x_{28}}-x_{12}\ptl_{x_{31}}-x_{14}\ptl_{x_{34}}
\\ &&+x_{23}\ptl_{x_{43}}+x_{26}\ptl_{x_{45}}+x_{29}\ptl_{x_{47}}
+x_{32}\ptl_{x_{48}}+x_{35}\ptl_{x_{50}}+x_{49}\ptl_{x_{56}},
\hspace{3.7cm}(2.85)
\end{eqnarray*}
\begin{eqnarray*}& &E_{(0,1,1,2,1,1)}|_V=x_2\ptl_{x_{10}}-x_4\ptl_{x_{14}}
+x_8\ptl_{x_{24}}-x_{11}\ptl_{x_{26}}+x_{16}\ptl_{x_{32}}-x_{18}\ptl_{x_{35}}
\\ &&-x_{22}\ptl_{x_{39}}+x_{25}\ptl_{x_{41}}-x_{31}\ptl_{x_{46}}
+x_{33}\ptl_{x_{49}}-x_{43}\ptl_{x_{53}}+x_{47}\ptl_{x_{55}},
\hspace{3.7cm}(2.86)
\end{eqnarray*}
\begin{eqnarray*}& &E_{(0,1,1,1,1,1,1)}|_V=x_1\ptl_{x_9}-x_5\ptl_{x_{17}}
+x_8\ptl_{x_{25}}+x_{13}\ptl_{x_{30}}+x_{16}\ptl_{x_{33}}+x_{19}\ptl_{x_{37}}
\\ &&+x_{20}\ptl_{x_{38}}+x_{24}\ptl_{x_{41}}+x_{27}\ptl_{x_{44}}
+x_{32}\ptl_{x_{49}}-x_{40}\ptl_{x_{52}}+x_{48}\ptl_{x_{56}},
\hspace{3.7cm}(2.87)
\end{eqnarray*}
\begin{eqnarray*}& &E_{(1,1,2,2,1)}|_V=-x_3\ptl_{x_{15}}+x_4\ptl_{x_{18}}
-x_5\ptl_{x_{20}}+x_7\ptl_{x_{23}}+x_{14}\ptl_{x_{35}}-x_{19}\ptl_{x_{40}}
\\ &&+x_{17}\ptl_{x_{38}}-x_{22}\ptl_{x_{43}}-x_{34}\ptl_{x_{50}}
+x_{37}\ptl_{x_{52}}-x_{39}\ptl_{x_{53}}+x_{42}\ptl_{x_{54}},
\hspace{3.7cm}(2.88)
\end{eqnarray*}
\begin{eqnarray*}& &E_{(1,1,1,2,1,1)}|_V=x_2\ptl_{x_{13}}-x_4\ptl_{x_{19}}
-x_6\ptl_{x_{24}}+x_9\ptl_{x_{26}}-x_{12}\ptl_{x_{32}}+x_{18}\ptl_{x_{40}}
\\ &&+x_{17}\ptl_{x_{39}}-x_{25}\ptl_{x_{45}}+x_{31}\ptl_{x_{48}}
-x_{33}\ptl_{x_{51}}-x_{38}\ptl_{x_{53}}+x_{44}\ptl_{x_{55}},
\hspace{3.7cm}(2.89)
\end{eqnarray*}
\begin{eqnarray*}& &E_{(1,1,1,1,1,1,1)}|_V=x_1\ptl_{x_{11}}-x_5\ptl_{x_{22}}
-x_6\ptl_{x_{25}}-x_{10}\ptl_{x_{30}}-x_{12}\ptl_{x_{33}}-x_{14}\ptl_{x_{37}}
\\ &&-x_{20}\ptl_{x_{43}}-x_{24}\ptl_{x_{45}}-x_{27}\ptl_{x_{47}}
-x_{32}\ptl_{x_{51}}-x_{35}\ptl_{x_{52}}+x_{46}\ptl_{x_{56}},
\hspace{3.7cm}(2.90)
\end{eqnarray*}
\begin{eqnarray*}& &E_{(0,1,1,2,2,1)}|_V=-x_2\ptl_{x_{12}}+x_3\ptl_{x_{14}}
-x_8\ptl_{x_{27}}+x_{11}\ptl_{x_{29}}+x_{13}\ptl_{x_{32}}-x_{15}\ptl_{x_{35}}
\\ &&+x_{22}\ptl_{x_{42}}-x_{25}\ptl_{x_{44}}-x_{28}\ptl_{x_{46}}
+x_{30}\ptl_{x_{49}}-x_{43}\ptl_{x_{54}}+x_{45}\ptl_{x_{55}},
\hspace{3.7cm}(2.91)
\end{eqnarray*}
\begin{eqnarray*}& &E_{(0,1,1,2,1,1,1)}|_V=-x_1\ptl_{x_{10}}+x_4\ptl_{x_{17}}
-x_8\ptl_{x_{28}}+x_{11}\ptl_{x_{30}}-x_{16}\ptl_{x_{36}}+x_{18}\ptl_{x_{38}}
\\ &&-x_{19}\ptl_{x_{39}}+x_{21}\ptl_{x_{41}}-x_{27}\ptl_{x_{46}}
+x_{29}\ptl_{x_{49}}-x_{40}\ptl_{x_{53}}+x_{47}\ptl_{x_{56}},
\hspace{3.7cm}(2.92)
\end{eqnarray*}
\begin{eqnarray*}& &E_{(1,1,2,2,1,1)}|_V=x_2\ptl_{x_{15}}-x_4\ptl_{x_{21}}
+x_5\ptl_{x_{24}}-x_7\ptl_{x_{26}}+x_{12}\ptl_{x_{35}}-x_{16}\ptl_{x_{40}}
\\ &&-x_{17}\ptl_{x_{41}}+x_{22}\ptl_{x_{45}}-x_{31}\ptl_{x_{50}}
+x_{33}\ptl_{x_{52}}-x_{36}\ptl_{x_{53}}+x_{42}\ptl_{x_{55}},
\hspace{3.7cm}(2.93)
\end{eqnarray*}
\begin{eqnarray*}& &E_{(1,1,1,2,2,1)}|_V=-x_2\ptl_{x_{16}}+x_3\ptl_{x_{19}}
+x_6\ptl_{x_{27}}-x_9\ptl_{x_{29}}-x_{10}\ptl_{x_{32}}+x_{15}\ptl_{x_{40}}
\\ &&-x_{17}\ptl_{x_{42}}+x_{25}\ptl_{x_{47}}+x_{28}\ptl_{x_{48}}
-x_{30}\ptl_{x_{51}}-x_{38}\ptl_{x_{54}}+x_{41}\ptl_{x_{55}},
\hspace{3.7cm}(2.94)
\end{eqnarray*}
\begin{eqnarray*}& &E_{(1,1,1,2,1,1,1)}|_V=-x_1\ptl_{x_{13}}+x_4\ptl_{x_{22}}
+x_6\ptl_{x_{28}}-x_9\ptl_{x_{30}}+x_{12}\ptl_{x_{36}}+x_{14}\ptl_{x_{39}}
\\ &&-x_{18}\ptl_{x_{43}}-x_{21}\ptl_{x_{45}}+x_{27}\ptl_{x_{48}}
-x_{29}\ptl_{x_{51}}-x_{35}\ptl_{x_{53}}+x_{44}\ptl_{x_{56}},
\hspace{3.7cm}(2.95)
\end{eqnarray*}
\begin{eqnarray*}& &E_{(0,1,1,2,2,1,1)}|_V=x_1\ptl_{x_{12}}-x_3\ptl_{x_{17}}
+x_8\ptl_{x_{31}}-x_{11}\ptl_{x_{33}}-x_{13}\ptl_{x_{36}}+x_{15}\ptl_{x_{38}}
\\ &&+x_{19}\ptl_{x_{42}}-x_{21}\ptl_{x_{44}}-x_{24}\ptl_{x_{46}}
+x_{26}\ptl_{x_{49}}-x_{40}\ptl_{x_{54}}+x_{45}\ptl_{x_{56}},
\hspace{3.7cm}(2.96)
\end{eqnarray*}
\begin{eqnarray*}& &E_{(1,1,2,2,2,1)}|_V=-x_2\ptl_{x_{18}}+x_3\ptl_{x_{21}}
-x_5\ptl_{x_{27}}+x_7\ptl_{x_{29}}+x_{10}\ptl_{x_{35}}-x_{13}\ptl_{x_{40}}
\\ &&+x_{17}\ptl_{x_{44}}-x_{22}\ptl_{x_{47}}-x_{28}\ptl_{x_{50}}
+x_{30}\ptl_{x_{52}}-x_{36}\ptl_{x_{54}}+x_{39}\ptl_{x_{55}},
\hspace{3.7cm}(2.97)
\end{eqnarray*}
\begin{eqnarray*}& &E_{(1,1,2,2,1,1,1)}|_V=-x_1\ptl_{x_{15}}+x_4\ptl_{x_{25}}
-x_5\ptl_{x_{28}}+x_7\ptl_{x_{30}}-x_{12}\ptl_{x_{38}}+x_{16}\ptl_{x_{43}}
\\ &&-x_{14}\ptl_{x_{41}}+x_{19}\ptl_{x_{45}}-x_{27}\ptl_{x_{50}}
+x_{29}\ptl_{x_{52}}-x_{32}\ptl_{x_{53}}+x_{42}\ptl_{x_{56}},
\hspace{3.7cm}(2.98)
\end{eqnarray*}
\begin{eqnarray*}& &E_{(1,1,1,2,2,1,1)}|_V=x_1\ptl_{x_{16}}-x_3\ptl_{x_{22}}
-x_6\ptl_{x_{31}}+x_9\ptl_{x_{33}}+x_{10}\ptl_{x_{36}}-x_{15}\ptl_{x_{43}}
\\ &&-x_{14}\ptl_{x_{42}}+x_{21}\ptl_{x_{47}}+x_{24}\ptl_{x_{48}}
-x_{26}\ptl_{x_{51}}-x_{35}\ptl_{x_{54}}+x_{41}\ptl_{x_{56}},
\hspace{3.7cm}(2.99)
\end{eqnarray*}
\begin{eqnarray*}& &E_{(0,1,1,2,2,2,1)}|_V=-x_1\ptl_{x_{14}}+x_2\ptl_{x_{17}}
-x_8\ptl_{x_{34}}+x_{11}\ptl_{x_{37}}+x_{13}\ptl_{x_{39}}-x_{15}\ptl_{x_{41}}
\\ &&+x_{16}\ptl_{x_{42}}-x_{18}\ptl_{x_{44}}-x_{20}\ptl_{x_{46}}
+x_{23}\ptl_{x_{49}}-x_{40}\ptl_{x_{55}}+x_{43}\ptl_{x_{56}},
\hspace{3.7cm}(2.100)
\end{eqnarray*}
\begin{eqnarray*}& &E_{(1,1,2,3,2,1)}|_V=x_2\ptl_{x_{20}}-x_3\ptl_{x_{24}}
+x_4\ptl_{x_{27}}-x_7\ptl_{x_{32}}+x_9\ptl_{x_{35}}-x_{11}\ptl_{x_{40}}
\\ &&-x_{17}\ptl_{x_{46}}+x_{22}\ptl_{x_{48}}-x_{25}\ptl_{x_{50}}
+x_{30}\ptl_{x_{53}}-x_{33}\ptl_{x_{54}}+x_{37}\ptl_{x_{55}},
\hspace{3.7cm}(2.101)
\end{eqnarray*}
\begin{eqnarray*}& &E_{(1,1,2,2,2,1,1)}|_V=x_1\ptl_{x_{18}}-x_3\ptl_{x_{25}}
+x_5\ptl_{x_{31}}-x_7\ptl_{x_{33}}-x_{10}\ptl_{x_{38}}+x_{13}\ptl_{x_{43}}
\\ &&+x_{14}\ptl_{x_{44}}-x_{19}\ptl_{x_{47}}-x_{24}\ptl_{x_{50}}
+x_{26}\ptl_{x_{52}}-x_{32}\ptl_{x_{54}}+x_{39}\ptl_{x_{56}},
\hspace{3.7cm}(2.102)
\end{eqnarray*}
\begin{eqnarray*}&
&E_{(1,1,1,2,2,2,1)}|_V=-x_1\ptl_{x_{19}}+x_2\ptl_{x_{22}}
+x_6\ptl_{x_{34}}-x_9\ptl_{x_{37}}-x_{10}\ptl_{x_{39}}+x_{15}\ptl_{x_{45}}
\\ &&-x_{12}\ptl_{x_{42}}+x_{18}\ptl_{x_{47}}+x_{20}\ptl_{x_{48}}
-x_{23}\ptl_{x_{51}}-x_{35}\ptl_{x_{55}}+x_{38}\ptl_{x_{56}},
\hspace{3.7cm}(2.103)
\end{eqnarray*}
\begin{eqnarray*}& &E_{(1,2,2,3,2,1)}|_V=-x_2\ptl_{x_{23}}+x_3\ptl_{x_{26}}
-x_4\ptl_{x_{29}}+x_5\ptl_{x_{32}}-x_6\ptl_{x_{35}}+x_8\ptl_{x_{40}}
\\ &&-x_{17}\ptl_{x_{49}}+x_{22}\ptl_{x_{51}}-x_{25}\ptl_{x_{52}}
+x_{28}\ptl_{x_{53}}-x_{31}\ptl_{x_{54}}+x_{34}\ptl_{x_{55}},
\hspace{3.7cm}(2.104)
\end{eqnarray*}
\begin{eqnarray*}& &E_{(1,1,2,3,2,1,1)}|_V=-x_1\ptl_{x_{20}}+x_3\ptl_{x_{28}}
-x_4\ptl_{x_{31}}+x_7\ptl_{x_{36}}-x_9\ptl_{x_{38}}+x_{11}\ptl_{x_{43}}
\\ &&-x_{14}\ptl_{x_{46}}+x_{19}\ptl_{x_{48}}-x_{21}\ptl_{x_{50}}
+x_{26}\ptl_{x_{53}}-x_{29}\ptl_{x_{54}}+x_{37}\ptl_{x_{56}},
\hspace{3.7cm}(2.105)
\end{eqnarray*}
\begin{eqnarray*}& &E_{(1,1,2,2,2,2,1)}|_V=-x_1\ptl_{x_{21}}+x_2\ptl_{x_{25}}
-x_5\ptl_{x_{34}}+x_7\ptl_{x_{37}}+x_{10}\ptl_{x_{41}}-x_{13}\ptl_{x_{45}}
\\ &&+x_{12}\ptl_{x_{44}}-x_{16}\ptl_{x_{47}}-x_{20}\ptl_{x_{50}}
+x_{23}\ptl_{x_{52}}-x_{32}\ptl_{x_{55}}+x_{36}\ptl_{x_{56}},
\hspace{3.7cm}(2.106)
\end{eqnarray*}
\begin{eqnarray*}& &E_{(1,2,2,3,2,1,1)}|_V=x_1\ptl_{x_{23}}-x_3\ptl_{x_{30}}
+x_4\ptl_{x_{33}}-x_5\ptl_{x_{36}}+x_6\ptl_{x_{38}}-x_8\ptl_{x_{43}}
\\ &&-x_{14}\ptl_{x_{49}}+x_{19}\ptl_{x_{51}}-x_{21}\ptl_{x_{52}}
+x_{24}\ptl_{x_{53}}-x_{27}\ptl_{x_{54}}+x_{34}\ptl_{x_{56}},
\hspace{3.7cm}(2.107)
\end{eqnarray*}
\begin{eqnarray*}& &E_{(1,1,2,3,2,2,1)}|_V=x_1\ptl_{x_{24}}-x_2\ptl_{x_{28}}
+x_4\ptl_{x_{34}}-x_7\ptl_{x_{39}}+x_9\ptl_{x_{41}}-x_{11}\ptl_{x_{45}}
\\ &&-x_{12}\ptl_{x_{46}}+x_{16}\ptl_{x_{48}}-x_{18}\ptl_{x_{50}}
+x_{23}\ptl_{x_{53}}-x_{29}\ptl_{x_{55}}+x_{33}\ptl_{x_{56}},
\hspace{3.7cm}(2.108)
\end{eqnarray*}
\begin{eqnarray*}& &E_{(1,2,2,3,2,2,1)}|_V=-x_1\ptl_{x_{26}}+x_2\ptl_{x_{30}}
-x_4\ptl_{x_{37}}+x_5\ptl_{x_{39}}-x_6\ptl_{x_{41}}+x_8\ptl_{x_{45}}
\\ &&-x_{12}\ptl_{x_{49}}+x_{16}\ptl_{x_{51}}-x_{18}\ptl_{x_{52}}
+x_{20}\ptl_{x_{53}}-x_{27}\ptl_{x_{55}}+x_{31}\ptl_{x_{56}},
\hspace{3.7cm}(2.109)
\end{eqnarray*}
\begin{eqnarray*}& &E_{(1,1,2,3,3,2,1)}|_V=-x_1\ptl_{x_{27}}+x_2\ptl_{x_{31}}
-x_3\ptl_{x_{34}}+x_7\ptl_{x_{42}}-x_9\ptl_{x_{44}}+x_{11}\ptl_{x_{47}}
\\ &&-x_{10}\ptl_{x_{46}}+x_{13}\ptl_{x_{48}}-x_{15}\ptl_{x_{50}}
+x_{23}\ptl_{x_{54}}-x_{26}\ptl_{x_{55}}+x_{30}\ptl_{x_{56}},
\hspace{3.6cm}(2.110)
\end{eqnarray*}
\begin{eqnarray*}&
&E_{(1,2,2,3,3,2,1)}|_V=x_1\ptl_{x_{29}}-x_2\ptl_{x_{33}}
+x_3\ptl_{x_{37}}-x_5\ptl_{x_{42}}+x_6\ptl_{x_{44}}-x_8\ptl_{x_{47}}
\\ &&-x_{10}\ptl_{x_{49}}+x_{13}\ptl_{x_{51}}-x_{15}\ptl_{x_{52}}
+x_{20}\ptl_{x_{54}}-x_{24}\ptl_{x_{55}}+x_{28}\ptl_{x_{56}},
\hspace{3.6cm}(2.111)
\end{eqnarray*}
\begin{eqnarray*}&
&E_{(1,2,2,4,3,2,1)}|_V=-x_1\ptl_{x_{32}}+x_2\ptl_{x_{36}}
-x_3\ptl_{x_{39}}+x_4\ptl_{x_{42}}-x_6\ptl_{x_{46}}+x_8\ptl_{x_{48}}
\\ &&-x_9\ptl_{x_{49}}+x_{11}\ptl_{x_{51}}-x_{15}\ptl_{x_{53}}
+x_{18}\ptl_{x_{54}}-x_{21}\ptl_{x_{55}}+x_{25}\ptl_{x_{56}},
\hspace{3.6cm}(2.112)
\end{eqnarray*}
\begin{eqnarray*}&
&E_{(1,2,3,4,3,2,1)}|_V=-x_1\ptl_{x_{35}}+x_2\ptl_{x_{38}}
-x_3\ptl_{x_{41}}+x_4\ptl_{x_{44}}-x_5\ptl_{x_{46}}+x_8\ptl_{x_{50}}
\\ &&-x_7\ptl_{x_{49}}+x_{11}\ptl_{x_{52}}-x_{13}\ptl_{x_{53}}
+x_{16}\ptl_{x_{54}}-x_{19}\ptl_{x_{55}}+x_{22}\ptl_{x_{56}},
\hspace{3.6cm}(2.113)
\end{eqnarray*}
\begin{eqnarray*}&
&E_{(2,2,3,4,3,2,1)}|_V=-x_1\ptl_{x_{40}}+x_2\ptl_{x_{43}}
-x_3\ptl_{x_{45}}+x_4\ptl_{x_{47}}-x_5\ptl_{x_{48}}+x_6\ptl_{x_{50}}
\\ &&-x_7\ptl_{x_{51}}+x_9\ptl_{x_{52}}-x_{10}\ptl_{x_{53}}
+x_{12}\ptl_{x_{54}}-x_{14}\ptl_{x_{55}}+x_{17}\ptl_{x_{56}}.
\hspace{3.8cm}(2.114)
\end{eqnarray*}

We define a symmetric linear operation $\tau$ on the space
$\sum_{i,j=1}^{56}\mbb{C}x_i\ptl_{x_j}$ by
$$\tau(x_i\ptl_{x_j})=x_j\ptl_{x_i}.\eqno(2.115)$$
 Then
$$E_{-\al}|_V=-\tau(E_{\al}|_V)\qquad\for\;\;\al\in\Phi_{E_7}^+\eqno(2.116)$$
by the second equations in (2.5) and (2.6). Moreover, we find
$$\al_r|_V=\sum_{i=1}^{28}a_{i,r}(x_i\ptl{x_i}-x_{\bar i}\ptl{x_{\bar
i}})\qquad\for\;\;r\in\ol{1,7},\eqno(2.117)$$ where $a_{i,r}$ are
constants given by the following table:
\begin{center}{\bf \large Table 1}\end{center}
\begin{center}\begin{tabular}{|r||r|r|r|r|r|r|r||r||r|r|r|r|r|r|r|}\hline
$i$&$a_{i,1}$&$a_{i,2}$&$a_{i,3}$&$a_{i,4}$&$a_{i,5}$&$a_{i,6}$&$a_{i,7}$&$i$&$a_{i,1}$&$a_{i,2}$&$a_{i,3}$&$a_{i,4}$&$a_{i,5}$&
$a_{i,6}$&$a_{i,7}$
\\\hline\hline 1&0&0&0&0&0&0&1&2&0&0&0&0&0&1&$-1$\\\hline
3&0&0&0&0&1&$-1$&0&4&$0$&0&0&1&$-1$&0&0
\\\hline 5&0&1&1&$-1$&0&0&0&6&1&1&$-1$&0&0&0&0 \\\hline
7&0&$-1$&$1$&0&0&0&0&8&$-1$&1&0&$0$&0&0&0\\\hline
9&1&$-1$&$-1$&1&0&0&0&10&1&0&0&$-1$&1&0&0\\\hline
11&$-1$&$-1$&0&1&$0$&0&0&12&1&0&0&0&$-1$&1&0\\\hline
13&$-1$&0&1&$-1$&1&$0$&0&14&1&0&0&0&0&$-1$&1\\\hline
15&0&0&$-1$&0&1&0&0&16&$-1$&0&1&0&$-1$&1&0\\\hline
17&1&0&0&0&0&0&$-1$&18&0&0&$-1$&1&$-1$&1&0\\\hline
19&$-1$&0&1&0&0&$-1$&1&20&0&1&0&$-1$&0&1&0\\\hline
21&0&0&$-1$&1&0&$-1$&1&22&$-1$&0&1&0&0&0&$-1$\\\hline
23&0&$-1$&0&0&0&1&0&24&0&1&0&$-1$&1&$-1$&1\\\hline
25&0&0&$-1$&1&0&0&$-1$&26&0&$-1$&0&0&1&$-1$&1\\\hline
27&0&1&0&0&$-1$&0&1&28&0&1
&0&$-1$&1&0&$-1$\\\hline\end{tabular}\end{center}

This completes the explicit construction of the basic irreducible
module of $E_7$ Lie algebra.

\section{Examples of Singular Vectors}

Now ${\cal A}=\mbb{C}[x_1,...,x_{56}]$ becomes a ${\cal
G}^{E_7}$-module via the differential operators in (2.52)-(2.117)

According to Table 1, we look for a singular vector of the form:
$$\zeta_1=c_1x_1x_{17}+c_2x_2x_{14}+c_3x_3x_{12}+c_4x_4x_{10}+c_5x_5x_9+c_6x_6x_7,\qquad c_i\in\mbb{C}.
\eqno(3.1)$$ Note that $E_{\al_1}(\zeta_1)=0$ naturally holds by
(2.52). Moreover, (2.53) gives
$$0=E_{\al_2}(\zeta_1)=(c_5+c_6)x_5x_6\lra c_6=-c_5,\eqno(3.2)$$
which also implies $E_{\al_3}(\zeta_1)=0$ by (2.54). Furthermore,
(2.55) implies
$$0=E_{\al_4}(\zeta_1)=(c_5-c_4)x_4x_9\lra c_5=c_4.\eqno(3.3)$$
Besides, (2.56) yields
$$0=E_{\al_5}(\zeta_1)=(c_4-c_3)x_3\ptl_{x_4}x_{10}\lra
c_4=c_3.\eqno(3.4)$$ According to (2.57), we have
$$0=E_{\al_6}(\zeta_1)=(c_3-c_2)x_2x_{12}\lra
c_3=c_2.\eqno(3.5)$$ In addition, (2.58) says that
$$0=E_{\al_7}(\zeta_1)=(c_2-c_1)x_1x_{14}\lra
c_2=c_1.\eqno(3.6)$$ Thus we have the singular vector
$$\zeta_1=x_1x_{17}+x_2x_{14}+x_3x_{12}+x_4x_{10}+x_5x_9-x_6x_7,
\eqno(3.7)$$ which generates an irreducible ${\cal G}^{E_7}$-module
$W$ isomorphic to the adjoint module.

Denote
$$\bar r=57-r\qquad\for\;\;r\in\ol{1,56}.\eqno(3.8)$$
Set
$$\mbb{A}={\cal
A}[\ptl_{x_1},...,\ptl_{x_{56}}].\eqno(3.9)$$ Define an associative
algebra isomorphism $\nu$ on $\mbb{A}$ by
$$\nu(x_r)=x_{\bar r},\;\nu(\ptl_{x_r})=\ptl_{x_{\bar
r}}\qquad\for\;\; r\in \ol{1,56}.\eqno(3.10)$$ According to
(2.52)-(2.58), (2.115) and (2.116), we define
\begin{eqnarray*}\td  E_{\al_1}|_V&=&\nu
(E_{\al_1}|_V)=E_{-\al_1}|_V= x_8\ptl_{x_6}+x_{11}\ptl_{x_9}
+x_{13}\ptl_{x_{10}}+x_{16}\ptl_{x_{12}}+x_{19}\ptl_{x_{14}}\\
&&+x_{22}\ptl_{x_{17}}
-x_{40}\ptl_{x_{35}}-x_{43}\ptl_{x_{38}}-x_{45}\ptl_{x_{41}}
-x_{47}\ptl_{x_{44}}-x_{48}\ptl_{x_{46}}-x_{51}\ptl_{x_{49}},\hspace{0.6cm}(3.11)
\end{eqnarray*}
\begin{eqnarray*}\td   E_{\al_2}|_V&=&\nu
(E_{\al_2}|_V)=-E_{-\al_2}|_V= x_7\ptl_{x_5}+x_9\ptl_{x_6}
+x_{11}\ptl_{x_8}-x_{23}\ptl_{x_{20}}-x_{26}\ptl_{x_{24}}\\
&&-x_{29}\ptl_{x_{27}}
-x_{30}\ptl_{x_{28}}-x_{33}\ptl_{x_{31}}-x_{37}\ptl_{x_{34}}
+x_{49}\ptl_{x_{46}}+x_{51}\ptl_{x_{48}}+x_{52}\ptl_{x_{50}},\hspace{0.6cm}(3.12)
\end{eqnarray*}
\begin{eqnarray*}\td   E_{\al_3}|_V&=&\nu
(E_{\al_3}|_V)=E_{-\al_3}|_V=x_6\ptl_{x_5}+x_9\ptl_{x_7}
+x_{15}\ptl_{x_{13}}+x_{18}\ptl_{x_{16}}+x_{21}\ptl_{x_{19}}\\
&&+x_{25}\ptl_{x_{22}}
-x_{35}\ptl_{x_{32}}-x_{38}\ptl_{x_{36}}-x_{41}\ptl_{x_{39}}
-x_{44}\ptl_{x_{42}}-x_{50}\ptl_{x_{48}}-x_{52}\ptl_{x_{51}},
\hspace{0.6cm}(3.13)
\end{eqnarray*}
\begin{eqnarray*}\td  E_{\al_4}|_V&=&\nu
(E_{\al_4}|_V)=-E_{-\al_4}|_V=x_5\ptl_{x_4}-x_{10}\ptl_{x_9}
-x_{13}\ptl_{x_{11}}-x_{20}\ptl_{x_{18}}-x_{24}\ptl_{x_{21}}\\
&&-x_{28}\ptl_{x_{25}}
-x_{32}\ptl_{x_{29}}-x_{36}\ptl_{x_{33}}-x_{39}\ptl_{x_{37}}
-x_{46}\ptl_{x_{44}}-x_{48}\ptl_{x_{47}}+x_{53}\ptl_{x_{52}},
\hspace{0.6cm}(3.14)
\end{eqnarray*}
\begin{eqnarray*}\td   E_{\al_5}|_V&=&\nu
(E_{\al_5}|_V)=-E_{-\al_5}|_V=x_4\ptl_{x_3}-x_{12}\ptl_{x_{10}}
-x_{16}\ptl_{x_{13}}-x_{18}\ptl_{x_{15}}-x_{27}\ptl_{x_{24}}\\
&&-x_{29}\ptl_{x_{26}}
-x_{31}\ptl_{x_{28}}-x_{33}\ptl_{x_{30}}-x_{42}\ptl_{x_{39}}
-x_{44}\ptl_{x_{41}}-x_{47}\ptl_{x_{45}}+x_{54}\ptl_{x_{53}},
\hspace{0.6cm}(3.15)
\end{eqnarray*}
\begin{eqnarray*}\td E_{\al_6}|_V&=&\nu
(E_{\al_6}|_V)=-E_{-\al_6}|_V=x_3\ptl_{x_2}-x_{14}\ptl_{x_{12}}
-x_{19}\ptl_{x_{16}}-x_{21}\ptl_{x_{18}}-x_{24}\ptl_{x_{20}}\\
&&-x_{26}\ptl_{x_{23}}
-x_{34}\ptl_{x_{31}}-x_{37}\ptl_{x_{33}}-x_{39}\ptl_{x_{36}}
-x_{41}\ptl_{x_{38}}-x_{45}\ptl_{x_{43}}+x_{55}\ptl_{x_{54}},
\hspace{0.6cm}(3.16)
\end{eqnarray*}
\begin{eqnarray*}\td E_{\al_7}|_V&=&\nu
(E_{\al_7}|_V)=-E_{-\al_7}|_V=x_2\ptl_{x_1}-x_{17}\ptl_{x_{14}}
-x_{22}\ptl_{x_{19}}-x_{25}\ptl_{x_{21}}-x_{28}\ptl_{x_{24}}\\
&&-x_{30}\ptl_{x_{26}}
-x_{31}\ptl_{x_{27}}-x_{33}\ptl_{x_{29}}-x_{36}\ptl_{x_{32}}
-x_{38}\ptl_{x_{35}}-x_{43}\ptl_{x_{40}}+x_{56}\ptl_{x_{55}}.
\hspace{0.6cm}(3.17)
\end{eqnarray*}

 Set
$$\zeta_2=\td
E_{\al_1}(\zeta_1)=x_1x_{22}+x_2x_{19}+x_3x_{16}+x_4x_{13}+x_5x_{11}-x_7x_8,
\eqno(3.18)$$
$$\zeta_3=\td
E_{\al_3}(\zeta_2)=x_1x_{25}+x_2x_{21}+x_3x_{18}+x_4x_{15}+x_6x_{11}-x_8x_9,
\eqno(3.19)$$
$$\zeta_4=\td
E_{\al_4}(\zeta_2)=-x_1x_{28}-x_2x_{24}-x_3x_{20}+x_5x_{15}-x_6x_{13}+x_8x_{10},
\eqno(3.20)$$
$$\zeta_5=\td
E_{\al_2}(\zeta_4)=x_1x_{30}+x_2x_{26}+x_3x_{23}+x_7x_{15}-x_9x_{13}+x_{10}x_{11},
\eqno(3.21)$$
$$\zeta_6=\td
E_{\al_5}(\zeta_4)=x_1x_{31}+x_2x_{27}-x_4x_{20}-x_5x_{18}+x_6x_{16}-x_8x_{12},
\eqno(3.22)$$
$$\zeta_7=\td
E_{\al_2}(\zeta_6)=-x_1x_{33}-x_2x_{29}+x_4x_{23}-x_7x_{18}+x_9x_{16}-x_{11}x_{12},
\eqno(3.23)$$
$$\zeta_8=\td
E_{\al_6}(\zeta_6)=-x_1x_{34}+x_3x_{27}+x_4x_{24}+x_5x_{21}-x_6x_{19}+x_8x_{14},
\eqno(3.24)$$
$$\zeta_9=\td
E_{\al_7}(\zeta_8)=-x_2x_{34}-x_3x_{31}-x_4x_{28}-x_5x_{25}+x_6x_{22}-x_8x_{17},
\eqno(3.25)$$
$$\zeta_{10}=\td
E_{\al_4}(\zeta_7)=x_1x_{36}+x_2x_{32}+x_5x_{23}+x_7x_{20}-x_{10}x_{16}+x_{12}x_{13},
\eqno(3.26)$$
$$\zeta_{11}=\td
E_{\al_2}(\zeta_8)=x_1x_{37}-x_3x_{29}-x_4x_{26}+x_7x_{21}-x_9x_{19}+x_{11}x_{14},
\eqno(3.27)$$
$$\zeta_{12}=\td
E_{\al_3}(\zeta_{10})=-x_1x_{38}-x_2x_{35}+x_6x_{23}+x_9x_{20}-x_{10}x_{18}+x_{12}x_{15},
\eqno(3.28)$$
$$\zeta_{13}=\td
E_{\al_2}(\zeta_9)=x_2x_{37}+x_3x_{33}+x_4x_{30}-x_7x_{25}+x_9x_{22}-x_{11}x_{17},
\eqno(3.29)$$
$$\zeta_{14}=\td
E_{\al_4}(\zeta_{11})=-x_1x_{39}+x_3x_{32}-x_5x_{26}-x_7x_{24}+x_{10}x_{19}-x_{13}x_{14},
\eqno(3.30)$$
$$\zeta_{15}=\td
E_{\al_4}(\zeta_{13})=-x_2x_{39}-x_3x_{36}+x_5x_{30}+x_7x_{28}-x_{10}x_{22}+x_{13}x_{17},
\eqno(3.31)$$
$$\zeta_{16}=\td
E_{\al_3}(\zeta_{14})=x_1x_{41}-x_3x_{35}-x_6x_{26}-x_9x_{24}+x_{10}x_{21}-x_{14}x_{15},
\eqno(3.32)$$
$$\zeta_{17}=\td
E_{\al_5}(\zeta_{14})=x_1x_{42}+x_4x_{32}+x_5x_{29}+x_7x_{27}-x_{12}x_{19}+x_{14}x_{16},
\eqno(3.33)$$
$$\zeta_{18}=\td
E_{\al_3}(\zeta_{15})=x_2x_{41}+x_3x_{38}+x_6x_{30}+x_9x_{28}-x_{10}x_{25}+x_{15}x_{17},
\eqno(3.34)$$
$$\zeta_{19}=\td
E_{\al_1}(\zeta_{12})=x_1x_{43}+x_2x_{40}+x_8x_{23}+x_{11}x_{20}-x_{13}x_{18}+x_{15}x_{16},
\eqno(3.35)$$
$$\zeta_{20}=\td
E_{\al_7}(\zeta_{17})=x_2x_{42}-x_4x_{36}-x_5x_{33}-x_7x_{31}+x_{12}x_{22}-x_{16}x_{17},
\eqno(3.36)$$
$$\zeta_{21}=\td
E_{\al_3}(\zeta_{17})=-x_1x_{44}-x_4x_{35}+x_6x_{29}+x_9x_{27}-x_{12}x_{21}+x_{14}x_{18},
\eqno(3.37)$$
$$\zeta_{22}=\td
E_{\al_6}(\zeta_{20})=x_3x_{42}+x_4x_{39}+x_5x_{37}+x_7x_{34}-x_{14}x_{22}+x_{17}x_{19},
\eqno(3.38)$$
$$\zeta_{38}=\td
E_{\al_1}(\zeta_{16})=-x_1x_{45}+x_3x_{40}-x_8x_{26}-x_{11}x_{24}+x_{13}x_{21}-x_{15}x_{19},
\eqno(3.39)$$
$$\zeta_{24}=\td
E_{\al_5}(\zeta_{18})=-x_2x_{44}+x_4x_{38}-x_6x_{33}-x_9x_{31}+x_{12}x_{25}-x_{17}x_{18},
\eqno(3.40)$$
$$\zeta_{25}=\td
E_{\al_4}(\zeta_{21})=x_1x_{46}-x_5x_{35}-x_6x_{32}-x_{10}x_{27}+x_{12}x_{24}-x_{14}x_{20},
\eqno(3.41)$$
$$\zeta_{26}=\td
E_{\al_1}(\zeta_{18})=-x_2x_{45}-x_3x_{43}+x_8x_{30}+x_{11}x_{28}-x_{13}x_{25}+x_{15}x_{22},
\eqno(3.42)$$
$$\zeta_{27}=\td
E_{\al_3}(\zeta_{22})=-x_3x_{44}-x_4x_{41}+x_6x_{37}+x_9x_{34}-x_{14}x_{25}+x_{17}x_{21},
\eqno(3.43)$$
$$\zeta_{28}=\td
E_{\al_1}(\zeta_{21})=x_1x_{47}+x_4x_{40}+x_8x_{29}+x_{11}x_{27}-x_{16}x_{21}+x_{18}x_{19},
\eqno(3.44)$$
$$\zeta_{29}=\td
E_{\al_4}(\zeta_{24})=x_2x_{46}+x_5x_{38}+x_6x_{36}+x_{10}x_{31}-x_{12}x_{28}+x_{17}x_{20},
\eqno(3.45)$$
$$\zeta_{30}=\td
E_{\al_1}(\zeta_{25})=-x_1x_{48}+x_5x_{40}-x_8x_{32}-x_{13}x_{27}+x_{16}x_{24}-x_{19}x_{20},
\eqno(3.46)$$
$$\zeta_{31}=\td
E_{\al_5}(\zeta_{26})=x_2x_{47}-x_4x_{43}-x_8x_{33}-x_{11}x_{31}+x_{16}x_{25}-x_{18}x_{22},
\eqno(3.47)$$
$$\zeta_{32}=\td
E_{\al_4}(\zeta_{27})=x_3x_{46}-x_5x_{41}-x_6x_{39}-x_{10}x_{34}+x_{14}x_{28}-x_{17}x_{24},
\eqno(3.48)$$
$$\zeta_{33}=\td
E_{\al_2}(\zeta_{25})=x_1x_{49}-x_7x_{35}-x_9x_{32}+x_{10}x_{29}-x_{12}x_{26}+x_{14}x_{23},
\eqno(3.49)$$
$$\zeta_{34}=\td
E_{\al_1}(\zeta_{29})=-x_2x_{48}-x_5x_{43}+x_8x_{36}+x_{13}x_{31}-x_{16}x_{28}+x_{20}x_{22},
\eqno(3.50)$$
$$\zeta_{35}=\td
E_{\al_1}(\zeta_{27})=x_3x_{47}+x_4x_{45}+x_8x_{37}+x_{11}x_{34}-x_{19}x_{25}+x_{21}x_{22},
\eqno(3.51)$$
$$\zeta_{36}=\td
E_{\al_5}(\zeta_{32})=x_4x_{46}+x_5x_{44}+x_6x_{42}+x_{12}x_{34}-x_{14}x_{31}+x_{17}x_{27},
\eqno(3.52)$$
$$\zeta_{37}=\td
E_{\al_3}(\zeta_{30})=x_1x_{50}+x_6x_{40}+x_8x_{35}-x_{15}x_{27}+x_{18}x_{24}-x_{20}x_{21},
\eqno(3.53)$$
$$\zeta_{38}=\td
E_{\al_2}(\zeta_{29})=x_2x_{49}+x_7x_{38}+x_9x_{36}-x_{10}x_{33}+x_{12}x_{30}-x_{17}x_{23},
\eqno(3.54)$$
$$\zeta_{39}=\td
E_{\al_1}(\zeta_{32})=-x_3x_{48}+x_5x_{45}-x_8x_{39}-x_{13}x_{34}+x_{19}x_{28}-x_{22}x_{24},
\eqno(3.55)$$
$$\zeta_{40}=\td
E_{\al_1}(\zeta_{33})=-x_1x_{51}+x_7x_{40}-x_{11}x_{32}+x_{13}x_{29}-x_{16}x_{26}+x_{19}x_{23},
\eqno(3.56)$$
$$\zeta_{41}=\td
E_{\al_1}(\zeta_{34})=x_2x_{50}-x_6x_{43}-x_8x_{38}+x_{15}x_{31}-x_{18}x_{28}+x_{20}x_{25},
\eqno(3.57)$$
$$\zeta_{42}=\td
E_{\al_2}(\zeta_{32})=x_3x_{49}-x_7x_{41}-x_9x_{39}+x_{10}x_{37}-x_{14}x_{30}+x_{17}x_{26},
\eqno(3.58)$$
$$\zeta_{43}=\td
E_{\al_1}(\zeta_{36})=-x_4x_{48}-x_5x_{47}+x_8x_{42}+x_{16}x_{34}-x_{19}x_{31}+x_{22}x_{27},
\eqno(3.59)$$
$$\zeta_{44}=\td
E_{\al_2}(\zeta_{37})=x_1x_{52}+x_9x_{40}+x_{11}x_{35}+x_{15}x_{29}-x_{18}x_{26}+x_{21}x_{23},
\eqno(3.60)$$
$$\zeta_{45}=\td
E_{\al_1}(\zeta_{38})=-x_2x_{51}-x_7x_{43}+x_{11}x_{36}-x_{13}x_{33}+x_{16}x_{30}-x_{22}x_{23},
\eqno(3.61)$$
$$\zeta_{46}=\td
E_{\al_3}(\zeta_{39})=x_3x_{50}+x_6x_{45}+x_8x_{41}-x_{15}x_{34}+x_{21}x_{28}-x_{24}x_{25},
\eqno(3.62)$$
$$\zeta_{47}=\td
E_{\al_2}(\zeta_{36})=x_4x_{49}+x_7x_{44}+x_9x_{42}-x_{12}x_{37}+x_{14}x_{33}-x_{17}x_{29},
\eqno(3.63)$$
$$\zeta_{48}=\td
E_{\al_4}(\zeta_{44})=x_1x_{53}-x_{10}x_{40}-x_{13}x_{35}-x_{15}x_{32}+x_{20}x_{26}-x_{23}x_{24},
\eqno(3.64)$$
$$\zeta_{49}=\td
E_{\al_2}(\zeta_{41})=x_2x_{52}-x_9x_{43}-x_{11}x_{38}-x_{15}x_{33}+x_{18}x_{30}-x_{23}x_{25},
\eqno(3.65)$$
$$\zeta_{50}=\td
E_{\al_1}(\zeta_{42})=-x_3x_{51}+x_7x_{45}-x_{11}x_{39}+x_{13}x_{37}-x_{19}x_{30}+x_{22}x_{26},
\eqno(3.66)$$
$$\zeta_{51}=\td
E_{\al_3}(\zeta_{43})=x_4x_{50}-x_6x_{47}-x_8x_{44}+x_{18}x_{34}-x_{21}x_{31}+x_{25}x_{27},
\eqno(3.67)$$
$$\zeta_{52}=\td
E_{\al_4}(\zeta_{47})=x_5x_{49}-x_7x_{46}-x_{10}x_{42}+x_{12}x_{39}-x_{14}x_{36}+x_{17}x_{32},
\eqno(3.68)$$
$$\zeta_{53}=\td
E_{\al_5}(\zeta_{48})=x_1x_{54}+x_{12}x_{40}+x_{16}x_{35}+x_{18}x_{32}-x_{20}x_{29}+x_{23}x_{27},
\eqno(3.69)$$
$$\zeta_{54}=\td
E_{\al_4}(\zeta_{49})=x_2x_{53}+x_{10}x_{43}+x_{13}x_{38}+x_{15}x_{36}-x_{20}x_{30}+x_{23}x_{28},
\eqno(3.70)$$
$$\zeta_{55}=\td
E_{\al_3}(\zeta_{50})=x_3x_{52}+x_9x_{45}+x_{11}x_{41}+x_{15}x_{37}-x_{21}x_{30}+x_{25}x_{26},
\eqno(3.71)$$
$$\zeta_{56}=\td
E_{\al_1}(\zeta_{47})=-x_4x_{51}-x_7x_{47}+x_{11}x_{42}-x_{16}x_{37}+x_{19}x_{33}-x_{22}x_{29},
\eqno(3.72)$$
$$\zeta_{57}=\td
E_{\al_4}(\zeta_{51})=x_5x_{50}+x_6x_{48}+x_8x_{46}-x_{20}x_{34}+x_{24}x_{31}-x_{27}x_{28},
\eqno(3.73)$$
$$\zeta_{58}=\td
E_{\al_3}(\zeta_{52})=x_6x_{49}-x_9x_{46}+x_{10}x_{44}-x_{12}x_{41}+x_{14}x_{38}-x_{17}x_{35},
\eqno(3.74)$$
$$\zeta_{59}=\td
E_{\al_6}(\zeta_{53})=x_1x_{55}-x_{14}x_{40}-x_{19}x_{35}-x_{21}x_{32}+x_{24}x_{29}-x_{26}x_{27},
\eqno(3.75)$$
$$\zeta_{60}=\td
E_{\al_5}(\zeta_{54})=x_2x_{54}-x_{12}x_{43}-x_{16}x_{38}-x_{18}x_{36}+x_{20}x_{33}-x_{23}x_{31},
\eqno(3.76)$$
$$\zeta_{61}=\td
E_{\al_4}(\zeta_{55})=x_3x_{53}-x_{10}x_{45}-x_{13}x_{41}-x_{15}x_{39}+x_{24}x_{30}-x_{26}x_{28},
\eqno(3.77)$$
$$\zeta_{62}=\td
E_{\al_3}(\zeta_{56})=x_4x_{52}-x_9x_{47}-x_{11}x_{44}-x_{18}x_{37}+x_{21}x_{33}-x_{25}x_{29},
\eqno(3.78)$$
$$\zeta_{63}=\td
E_{\al_1}(\zeta_{52})=-x_5x_{51}+x_7x_{48}-x_{13}x_{42}+x_{16}x_{39}-x_{19}x_{36}+x_{22}x_{32},
\eqno(3.79)$$
\begin{eqnarray*}\zeta_{64}&=&\td
E_{\al_7}(\zeta_{59})=x_1x_{56}+x_2x_{55}+x_{14}x_{43}+x_{17}x_{40}+x_{19}x_{38}
\\ & &+x_{22}x_{35}+x_{21}x_{36}+x_{25}x_{32}-x_{24}x_{33}-x_{28}x_{29}+x_{26}x_{31}+x_{27}x_{30},
\hspace{2.1cm}(3.80)
\end{eqnarray*}
\begin{eqnarray*}\zeta_{65}&=&\td
E_{\al_6}(\zeta_{60})=x_2x_{55}+x_3x_{54}+x_{12}x_{45}+x_{14}x_{43}+x_{16}x_{41}+x_{19}x_{38}
\\ & &+x_{18}x_{39}+x_{21}x_{36}-x_{20}x_{37}-x_{24}x_{33}+x_{23}x_{34}+x_{26}x_{31},\hspace{3.7cm}(3.81)\end{eqnarray*}
\begin{eqnarray*}\zeta_{66}=\td
E_{\al_5}(\zeta_{61})&=&x_3x_{54}+x_4x_{53}+x_{10}x_{47}+x_{12}x_{45}+x_{13}x_{44}
+x_{15}x_{42}\\&&+x_{16}x_{41}+x_{18}x_{39}-x_{24}x_{33}-x_{27}x_{30}+x_{26}x_{31}+x_{28}x_{29},
\hspace{1.6cm}(3.82)\end{eqnarray*}
\begin{eqnarray*}\zeta_{67}=\td
E_{\al_4}(\zeta_{62})&=&x_4x_{53}+x_5x_{52}+x_9x_{48}+x_{10}x_{47}+x_{11}x_{46}+x_{13}x_{44}\\
&&+x_{18}x_{39}+x_{20}x_{37}-x_{21}x_{36}-x_{24}x_{33}+x_{25}x_{32}+x_{28}x_{29},
\hspace{1.7cm}(3.83)\end{eqnarray*}
\begin{eqnarray*}\zeta_{68}=\td
E_{\al_2}(\zeta_{57})&=&x_5x_{52}+x_6x_{51}+x_7x_{50}+x_8x_{49}+x_9x_{48}+x_{11}x_{46}\\
&&+x_{20}x_{37}+x_{23}x_{34}-x_{24}x_{33}-x_{26}x_{31}+x_{27}x_{30}+x_{28}x_{29},\hspace{1.7cm}(3.84)\end{eqnarray*}
\begin{eqnarray*}\zeta_{69}=\td
E_{\al_3}(\zeta_{63})&=&x_5x_{52}-x_6x_{51}-x_7x_{50}+x_9x_{48}+x_{13}x_{44}-x_{15}x_{42}\\
& &-x_{16}x_{41}+x_{18}x_{39}+x_{19}x_{38}-x_{21}x_{36}-x_{22}x_{35}
 +x_{25}x_{32},\hspace{1.7cm}(3.85)\end{eqnarray*}
\begin{eqnarray*}\zeta_{70}=\td
E_{\al_1}(\zeta_{58})&=&-x_6x_{51}+x_8x_{49}+x_9x_{48}-x_{10}x_{47}-x_{11}x_{46}+x_{13}x_{44}
\\ & &+x_{12}x_{45}-x_{14}x_{43}-x_{16}x_{41}+x_{17}x_{40}+x_{19}x_{38}-x_{22}x_{35}.\hspace{1.8cm}(3.86)\end{eqnarray*}

By calculating dimensions, we find that the subspace of homogeneous
quadratic polynomials in ${\cal A}$ is a direct sum of  irreducible
${\cal G}^{E_7}$-submodules generated by $x_1^2$ of weight $2\lmd_7$
and by $\zeta_1$ of weight $\lmd_1$. Since $\nu(\zeta_1)$ is a
lowest-weight vector of weight $-\lmd_1$, we must have
$\nu(\zeta_1)\in W$. Set
$$\zeta_i=\tau(\zeta_{134-i})\qquad\for\;\;i\in\ol{71,133}\eqno(3.87)$$
(cf. (3.8)-(3.10)). By (3.11)-(3.79) and induction on $134-i$, we
can prove $\zeta_i\in W$. Comparing weights and dimensions, we get
that $\{\zeta_r\mid r\in\ol{1,133}\}$ forms a basis of $W$.
Moreover,
$$\al_r|_W=\sum_{i=1}^{63}b_{i,r}(\zeta_i\ptl_{\zeta_i}-\zeta_{134-i}\ptl_{\zeta_{134-i}})
\qquad\for\;\;r\in\ol{1,7},\eqno(3.88)$$ where $b_{i,r}$ are given
by the Table 2 below. Note that for $\al,\be\in \Phi_{E_7}$,
$$\al+\be\in \Phi_{E_7}\Leftrightarrow (\al,\be)=-1.\eqno(3.89)$$
Thus
$$E_{\al_r}(\zeta_i)\neq 0 \Leftrightarrow
b_{i,r}=-1\qquad\for\;\;r\in\ol{1,7},\;i\in\ol{1,63}.\eqno(3.90)$$
Set
$$W_1=\sum_{i=1}^{63}\mbb{C}\zeta_i.\eqno(3.91)$$
\newpage

\begin{center}{\bf \large Table 2}\end{center}
\begin{center}\begin{tabular}{|r||r|r|r|r|r|r|r||r||r|r|r|r|r|r|r|}\hline
$i$&$b_{i,1}$&$b_{i,2}$&$b_{i,3}$&$b_{i,4}$&$b_{i,5}$&$b_{i,6}$&$b_{i,7}$&$i$&$b_{i,1}$&$b_{i,2}$&$b_{i,3}$&$b_{i,4}$&$b_{i,5}$&
$b_{i,6}$&$b_{i,7}$
\\\hline\hline 1&1&0&0&0&0&0&0&2&$-1$&0&1&0&0&0&0\\\hline
3&0&0&$-1$&1&0&0&0&4&0&1 &0&$-1$&1&0&0
\\\hline 5&0&$-1$&0&0&1&0&0&6&0&1&0&0&$-1$&1&0 \\\hline
7&0&$-1$&0&1&$-1$&1&0&8&0&1&0&0&0&$-1$&1\\\hline
9&0&1&0&0&0&0&$-1$&10&0&0&1&$-1$&0&1&0\\\hline
11&0&$-1$&0&1&0&$-1$&1&12&1&0&$-1$&0&0&1&0\\\hline
13&0&$-1$&0&1&0&0&$-1$&14&0&0&1&$-1$&1&$-1$&1\\\hline
15&0&0&1&$-1$&1&0&$-1$&16&1&0&$-1$&0&1&$-1$&1\\\hline
17&0&0&1&0&$-1$&0&1&18&1&0&$-1$&0&1&0&$-1$\\\hline
19&$-1$&0&0&0&0&1&0&20&0&0&1&0&$-1$&1&$-1$\\\hline
21&1&0&$-1$&1&$-1$&0&1&22&0&0&1&0&0&$-1$&0\\\hline
23&$-1$&0&0&0&1&$-1$&1&24&1&0&$-1$&1&$-1$&1&$-1$\\\hline
25&1&1&0&$-1$&0&0&1&26&$-1$&0&0&0&1&0&$-1$\\\hline
27&1&0&$-1$&1&0&$-1$&0&28&$-1$&0&0&1&$-1$&0&1
\\\hline 29&1&1&0&$-1$&0&1&$-1$&30&$-1$&1&1&$-1$&0&0&1\\\hline
31&$-1$&0&0&1&$-1$&1&$-1$&32&1&1&0&$-1$&1&$-1$&0
\\\hline 33&1&$-1$&0&$0$&0&0&1&34&$-1$&1&1&$-1$&0&1&$-1$ \\\hline
35&$-1$&0&0&1&0&$-1$&0&36&1&1&0&0&$-1$&0&0\\\hline
37&0&1&$-1$&0&0&0&1&38&1&$-1$&0&$0$&0&1&$-1$\\\hline
39&$-1$&1&1&$-1$&1&$-1$&0&40&$-1$&$-1$&1&0&0&0&1\\\hline
41&0&1&$-1$&0&0&1&$-1$&42&1&$-1$&0&$0$&1&$-1$&0\\\hline
43&$-1$&1&1&0&$-1$&0&0&44&0&$-1$&$-1$&1&0&0&1\\\hline
45&$-1$&$-1$&1&0&0&1&$-1$&46&0&1&$-1$&0&1&$-1$&0\\\hline
47&1&$-1$&0&1&$-1$&0&0&48&0&0&0&$-1$&1&0&1\\\hline
49&0&$-1$&$-1$&1&0&1&$-1$&50&$-1$&$-1$&1&0&1&$-1$&0\\\hline
51&0&1&$-1$&1&$-1$&0&0&52&1&0&1&$-1$&0&0&0\\\hline
53&0&0&0&0&$-1$&1&1&54&$0$&0&0&$-1$&1&1&$-1$\\\hline
55&0&$-1$&$-1$&1&1&$-1$&0&56&$-1$&$-1$&1&1&$-1$&0&0
\\\hline 57&0&2&0&$-1$&0&0&0&58&2&0&$-1$&0&0&0&0\\\hline
59&0&0&0&0&0&$-1$&2&60&0&0&0&0&$-1$&2&$-1$\\\hline
61&0&0&0&$-1$&2&$-1$&0&62&0&$-1$&$-1$&2&$-1$&0&0\\\hline
63&$-1$&0&2&$-1$&0&0&0&&&&&&&&\\\hline
\end{tabular}\end{center}

 Based on the above fact and Table 2, we find
\begin{eqnarray*}E_{\al_1}|_{W_1}&=&-\zeta_1\ptl_{\zeta_2}-\zeta_{12}\ptl_{\zeta_{19}}
-\zeta_{16}\ptl_{\zeta_{23}}-\zeta_{18}\ptl_{\zeta_{26}}-\zeta_{21}\ptl_{\zeta_{28}}-\zeta_{25}\ptl_{\zeta_{30}}
\\&&-\zeta_{24}\ptl_{\zeta_{31}}-\zeta_{29}\ptl_{\zeta_{34}}-\zeta_{27}\ptl_{\zeta_{35}}
-\zeta_{32}\ptl_{\zeta_{39}}-\zeta_{33}\ptl_{\zeta_{40}}
\\ &&-\zeta_{36}\ptl_{\zeta_{43}}-\zeta_{38}\ptl_{\zeta_{45}}-\zeta_{42}\ptl_{\zeta_{50}}
-\zeta_{47}\ptl_{\zeta_{56}}-\zeta_{52}\ptl_{\zeta_{63}},\hspace{4.4cm}(3.92)
\end{eqnarray*}
\begin{eqnarray*}E_{\al_2}|_{W_1}&=&\zeta_4\ptl_{\zeta_5}+\zeta_6\ptl_{\zeta_7}
+\zeta_8\ptl_{\zeta_{11}}+\zeta_9\ptl_{\zeta_{13}}+\zeta_{25}\ptl_{\zeta_{33}}+\zeta_{29}\ptl_{\zeta_{38}}
\\&&+\zeta_{30}\ptl_{\zeta_{40}}+\zeta_{32}\ptl_{\zeta_{42}}+\zeta_{37}\ptl_{\zeta_{44}}
+\zeta_{34}\ptl_{\zeta_{45}}+\zeta_{36}\ptl_{\zeta_{47}}
\\ &&+\zeta_{41}\ptl_{\zeta_{49}}+\zeta_{39}\ptl_{\zeta_{50}}+\zeta_{46}\ptl_{\zeta_{55}}
+\zeta_{43}\ptl_{\zeta_{56}}+\zeta_{51}\ptl_{\zeta_{62}},\hspace{4.4cm}(3.93)
\end{eqnarray*}
\begin{eqnarray*}E_{\al_3}|_{W_1}&=&-\zeta_2\ptl_{\zeta_3}-\zeta_{10}\ptl_{\zeta_{12}}
-\zeta_{14}\ptl_{\zeta_{16}}-\zeta_{15}\ptl_{\zeta_{18}}-\zeta_{17}\ptl_{\zeta_{21}}-\zeta_{20}\ptl_{\zeta_{24}}
\\&&-\zeta_{22}\ptl_{\zeta_{27}}-\zeta_{30}\ptl_{\zeta_{37}}-\zeta_{34}\ptl_{\zeta_{41}}
-\zeta_{40}\ptl_{\zeta_{44}}-\zeta_{39}\ptl_{\zeta_{46}}
\\ &&-\zeta_{45}\ptl_{\zeta_{49}}-\zeta_{43}\ptl_{\zeta_{51}}-\zeta_{50}\ptl_{\zeta_{55}}
-\zeta_{52}\ptl_{\zeta_{58}}-\zeta_{56}\ptl_{\zeta_{62}},\hspace{4.4cm}(3.94)
\end{eqnarray*}
\begin{eqnarray*}E_{\al_4}|_{W_1}&=&\zeta_3\ptl_{\zeta_4}+\zeta_7\ptl_{\zeta_{10}}
+\zeta_{11}\ptl_{\zeta_{14}}+\zeta_{13}\ptl_{\zeta_{15}}+\zeta_{21}\ptl_{\zeta_{25}}+\zeta_{24}\ptl_{\zeta_{29}}
\\&&+\zeta_{28}\ptl_{\zeta_{30}}+\zeta_{27}\ptl_{\zeta_{32}}+\zeta_{31}\ptl_{\zeta_{34}}
+\zeta_{35}\ptl_{\zeta_{39}}+\zeta_{44}\ptl_{\zeta_{48}}
\\ &&+\zeta_{47}\ptl_{\zeta_{52}}+\zeta_{49}\ptl_{\zeta_{54}}+\zeta_{51}\ptl_{\zeta_{57}}
+\zeta_{55}\ptl_{\zeta_{61}}+\zeta_{56}\ptl_{\zeta_{63}},\hspace{4.4cm}(3.95)
\end{eqnarray*}
\begin{eqnarray*}E_{\al_5}|_{W_1}&=&\zeta_4\ptl_{\zeta_6}+\zeta_5\ptl_{\zeta_7}
+\zeta_{14}\ptl_{\zeta_{17}}+\zeta_{15}\ptl_{\zeta_{20}}+\zeta_{16}\ptl_{\zeta_{21}}+\zeta_{18}\ptl_{\zeta_{24}}
\\&&+\zeta_{23}\ptl_{\zeta_{28}}+\zeta_{26}\ptl_{\zeta_{31}}+\zeta_{32}\ptl_{\zeta_{36}}
+\zeta_{39}\ptl_{\zeta_{43}}+\zeta_{42}\ptl_{\zeta_{47}}
\\ &&+\zeta_{46}\ptl_{\zeta_{51}}+\zeta_{48}\ptl_{\zeta_{53}}+\zeta_{50}\ptl_{\zeta_{56}}
+\zeta_{54}\ptl_{\zeta_{60}}+\zeta_{55}\ptl_{\zeta_{62}},\hspace{4.4cm}(3.96)
\end{eqnarray*}
\begin{eqnarray*}E_{\al_6}|_{W_1}&=&\zeta_6\ptl_{\zeta_8}+\zeta_7\ptl_{\zeta_{11}}
+\zeta_{10}\ptl_{\zeta_{14}}+\zeta_{12}\ptl_{\zeta_{16}}+\zeta_{20}\ptl_{\zeta_{22}}+\zeta_{19}\ptl_{\zeta_{23}}
\\&&+\zeta_{24}\ptl_{\zeta_{27}}+\zeta_{29}\ptl_{\zeta_{32}}+\zeta_{31}\ptl_{\zeta_{35}}
+\zeta_{34}\ptl_{\zeta_{39}}+\zeta_{38}\ptl_{\zeta_{42}}
\\ &&+\zeta_{41}\ptl_{\zeta_{46}}+\zeta_{45}\ptl_{\zeta_{50}}+\zeta_{49}\ptl_{\zeta_{55}}
+\zeta_{53}\ptl_{\zeta_{59}}+\zeta_{54}\ptl_{\zeta_{61}},\hspace{4.4cm}(3.97)
\end{eqnarray*}
\begin{eqnarray*}E_{\al_7}|_{W_1}&=&\zeta_8\ptl_{\zeta_9}+\zeta_{11}\ptl_{\zeta_{13}}
+\zeta_{14}\ptl_{\zeta_{15}}+\zeta_{16}\ptl_{\zeta_{18}}+\zeta_{17}\ptl_{\zeta_{20}}+\zeta_{21}\ptl_{\zeta_{24}}
\\&&+\zeta_{23}\ptl_{\zeta_{26}}+\zeta_{25}\ptl_{\zeta_{29}}+\zeta_{28}\ptl_{\zeta_{31}}
+\zeta_{30}\ptl_{\zeta_{34}}+\zeta_{33}\ptl_{\zeta_{38}}
\\ &&+\zeta_{37}\ptl_{\zeta_{41}}+\zeta_{40}\ptl_{\zeta_{45}}+\zeta_{44}\ptl_{\zeta_{49}}
+\zeta_{48}\ptl_{\zeta_{54}}+\zeta_{53}\ptl_{\zeta_{60}}.\hspace{4.4cm}(3.98)
\end{eqnarray*}

Similarly we have:
\begin{eqnarray*}\td  E_{\al_1}|_{W_1}&=&\zeta_2\ptl_{\zeta_1}+\zeta_{19}\ptl_{\zeta_{12}}
+\zeta_{23}\ptl_{\zeta_{16}}+\zeta_{26}\ptl_{\zeta_{18}}+\zeta_{28}\ptl_{\zeta_{21}}+\zeta_{30}\ptl_{\zeta_{25}}
\\&&+\zeta_{31}\ptl_{\zeta_{24}}+\zeta_{34}\ptl_{\zeta_{29}}+\zeta_{35}\ptl_{\zeta_{27}}
+\zeta_{39}\ptl_{\zeta_{32}}+\zeta_{40}\ptl_{\zeta_{33}}
+\zeta_{43}\ptl_{\zeta_{36}}\\
&&+\zeta_{45}\ptl_{\zeta_{38}}+\zeta_{50}\ptl_{\zeta_{42}}
+\zeta_{56}\ptl_{\zeta_{47}}+\zeta_{63}\ptl_{\zeta_{52}}+\zeta_{70}\ptl_{\zeta_{58}},\hspace{4.3cm}(3.99)
\end{eqnarray*}
\begin{eqnarray*}\td  E_{\al_2}|_{W_1}&=&\zeta_5\ptl_{\zeta_4}+\zeta_7\ptl_{\zeta_6}
+\zeta_{11}\ptl_{\zeta_8}+\zeta_{13}\ptl_{\zeta_9}+\zeta_{33}\ptl_{\zeta_{25}}+\zeta_{38}\ptl_{\zeta_{29}}
\\&&+\zeta_{40}\ptl_{\zeta_{30}}+\zeta_{42}\ptl_{\zeta_{32}}+\zeta_{44}\ptl_{\zeta_{37}}
+\zeta_{45}\ptl_{\zeta_{34}}+\zeta_{47}\ptl_{\zeta_{36}}
+\zeta_{49}\ptl_{\zeta_{41}}\\
&&+\zeta_{50}\ptl_{\zeta_{39}}+\zeta_{55}\ptl_{\zeta_{46}}
+\zeta_{56}\ptl_{\zeta_{43}}+\zeta_{62}\ptl_{\zeta_{51}}+\zeta_{68}\ptl_{\zeta_{57}},\hspace{4.2cm}(3.100)
\end{eqnarray*}
\begin{eqnarray*}\td  E_{\al_3}|_{W_1}&=&\zeta_3\ptl_{\zeta_2}+\zeta_{12}\ptl_{\zeta_{10}}
+\zeta_{16}\ptl_{\zeta_{14}}+\zeta_{18}\ptl_{\zeta_{15}}+\zeta_{21}\ptl_{\zeta_{17}}+\zeta_{24}\ptl_{\zeta_{20}}
\\&&+\zeta_{27}\ptl_{\zeta_{22}}+\zeta_{37}\ptl_{\zeta_{30}}+\zeta_{41}\ptl_{\zeta_{34}}
+\zeta_{44}\ptl_{\zeta_{40}}+\zeta_{46}\ptl_{\zeta_{39}}
+\zeta_{49}\ptl_{\zeta_{45}}\\
&&+\zeta_{51}\ptl_{\zeta_{43}}+\zeta_{55}\ptl_{\zeta_{50}}
+\zeta_{58}\ptl_{\zeta_{52}}+\zeta_{62}\ptl_{\zeta_{56}}+\zeta_{69}\ptl_{\zeta_{63}},\hspace{4.2cm}(3.101)
\end{eqnarray*}
\begin{eqnarray*}\td  E_{\al_4}|_{W_1}&=&\zeta_4\ptl_{\zeta_3}+\zeta_{10}\ptl_{\zeta_7}
+\zeta_{14}\ptl_{\zeta_{11}}+\zeta_{15}\ptl_{\zeta_{13}}+\zeta_{25}\ptl_{\zeta_{21}}+\zeta_{29}\ptl_{\zeta_{24}}
\\&&+\zeta_{30}\ptl_{\zeta_{28}}+\zeta_{32}\ptl_{\zeta_{27}}+\zeta_{34}\ptl_{\zeta_{31}}
+\zeta_{39}\ptl_{\zeta_{35}}+\zeta_{48}\ptl_{\zeta_{44}}
+\zeta_{52}\ptl_{\zeta_{47}}\\
&&+\zeta_{54}\ptl_{\zeta_{49}}+\zeta_{57}\ptl_{\zeta_{51}}
+\zeta_{61}\ptl_{\zeta_{55}}+\zeta_{63}\ptl_{\zeta_{56}}+\zeta_{67}\ptl_{\zeta_{62}},
\hspace{4.2cm}(3.102)
\end{eqnarray*}
\begin{eqnarray*}\td  E_{\al_5}|_{W_1}&=&\zeta_6\ptl_{\zeta_4}+\zeta_7\ptl_{\zeta_5}
+\zeta_{17}\ptl_{\zeta_{14}}+\zeta_{20}\ptl_{\zeta_{15}}+\zeta_{21}\ptl_{\zeta_{16}}+\zeta_{24}\ptl_{\zeta_{18}}
\\&&+\zeta_{28}\ptl_{\zeta_{23}}+\zeta_{31}\ptl_{\zeta_{26}}+\zeta_{36}\ptl_{\zeta_{32}}
+\zeta_{43}\ptl_{\zeta_{39}}+\zeta_{47}\ptl_{\zeta_{42}}
+\zeta_{51}\ptl_{\zeta_{46}}\\
&&+\zeta_{53}\ptl_{\zeta_{48}}+\zeta_{56}\ptl_{\zeta_{50}}
+\zeta_{60}\ptl_{\zeta_{54}}+\zeta_{62}\ptl_{\zeta_{55}}+\zeta_{66}\ptl_{\zeta_{61}},\hspace{4.2cm}(3.103)
\end{eqnarray*}
\begin{eqnarray*}\td  E_{\al_6}|_{W_1}&=&\zeta_8\ptl_{\zeta_6}+\zeta_{11}\ptl_{\zeta_7}
+\zeta_{14}\ptl_{\zeta_{10}}+\zeta_{16}\ptl_{\zeta_{14}}+\zeta_{22}\ptl_{\zeta_{20}}+\zeta_{23}\ptl_{\zeta_{19}}
\\&&+\zeta_{27}\ptl_{\zeta_{24}}+\zeta_{32}\ptl_{\zeta_{29}}+\zeta_{35}\ptl_{\zeta_{31}}
+\zeta_{39}\ptl_{\zeta_{34}}+\zeta_{42}\ptl_{\zeta_{38}}
+\zeta_{46}\ptl_{\zeta_{41}}\\
& &+\zeta_{50}\ptl_{\zeta_{45}}+\zeta_{55}\ptl_{\zeta_{49}}
+\zeta_{59}\ptl_{\zeta_{53}}+\zeta_{61}\ptl_{\zeta_{54}}+\zeta_{65}\ptl_{\zeta_{60}},\hspace{4.2cm}(3.104)
\end{eqnarray*}
\begin{eqnarray*}\td  E_{\al_7}|_{W_1}&=&\zeta_9\ptl_{\zeta_8}+\zeta_{13}\ptl_{\zeta_{11}}
+\zeta_{15}\ptl_{\zeta_{14}}+\zeta_{18}\ptl_{\zeta_{16}}+\zeta_{20}\ptl_{\zeta_{17}}+\zeta_{24}\ptl_{\zeta_{21}}
\\&&+\zeta_{26}\ptl_{\zeta_{23}}+\zeta_{29}\ptl_{\zeta_{25}}+\zeta_{31}\ptl_{\zeta_{28}}
+\zeta_{34}\ptl_{\zeta_{30}}+\zeta_{38}\ptl_{\zeta_{33}}
+\zeta_{41}\ptl_{\zeta_{37}}\\
&&+\zeta_{45}\ptl_{\zeta_{40}}+\zeta_{49}\ptl_{\zeta_{44}}
+\zeta_{54}\ptl_{\zeta_{48}}+\zeta_{60}\ptl_{\zeta_{53}}+\zeta_{64}\ptl_{\zeta_{59}}.\hspace{4.3cm}(3.105)
\end{eqnarray*}

We define a linear transformation $\Im$ on the space
$$\sum_{i,j=1}^{133}\mbb{C}\zeta_i\ptl_{\zeta_j}\eqno(3.106)$$
by
$$\Im(\zeta_i\ptl_{\zeta_j})=\zeta_{134-j}\ptl_{\zeta_{134-i}}.\eqno(3.107)$$
By symmetry, we have
$$E_{\al_1}|_W=(1-\Im)(E_{\al_1}|_{W_1})+\zeta_{70}\ptl_{\zeta_{76}}-\zeta_{58}(2\ptl_{\zeta_{70}}+\ptl_{\zeta_{69}}),
\eqno(3.108)$$
$$E_{\al_2}|_W=(1+\Im)(E_{\al_2}|_{W_1})+\zeta_{68}\ptl_{\zeta_{77}}+\zeta_{57}(2\ptl_{\zeta_{68}}+\ptl_{\zeta_{67}}),
\eqno(3.109)$$
$$E_{\al_3}|_W=(1-\Im)(E_{\al_3}|_{W_1})+\zeta_{69}\ptl_{\zeta_{71}}-
\zeta_{63}(\ptl_{\zeta_{70}}+2\ptl_{\zeta_{69}}+\ptl_{\zeta_{67}}),\eqno(3.110)$$
$$E_{\al_4}|_W=(1+\Im)(E_{\al_4}|_{W_1})+\zeta_{67}\ptl_{\zeta_{72}}+\zeta_{62}(\ptl_{\zeta_{69}}+\ptl_{\zeta_{68}}
+2\ptl_{\zeta_{67}}+\ptl_{\zeta_{66}}),\eqno(3.111)$$
$$E_{\al_5}|_W=(1+\Im)(E_{\al_5}|_{W_1})+\zeta_{66}\ptl_{\zeta_{73}}+\zeta_{61}(\ptl_{\zeta_{67}}
+2\ptl_{\zeta_{66}}+\ptl_{\zeta_{65}}),\eqno(3.112)$$
$$E_{\al_6}|_W=(1+\Im)(E_{\al_6}|_{W_1})+\zeta_{65}\ptl_{\zeta_{74}}+
\zeta_{60}(\ptl_{\zeta_{66}}+2\ptl_{\zeta_{65}}+\ptl_{\zeta_{64}}),
\eqno(3.113)$$
$$E_{\al_7}|_W=(1+\Im)(E_{\al_7}|_{W_1})+\zeta_{64}\ptl_{\zeta_{75}}+
\zeta_{59}(\ptl_{\zeta_{65}}+2\ptl_{\zeta_{64}}). \eqno(3.114)$$

 According Table 1 and Table 2, we look
for the singular vector of weight $\lmd_7$:
\begin{eqnarray*}& &\vt=a_1x_1\zeta_{64}+a_2x_2\zeta_{59}+a_3x_3\zeta_{53}+a_4x_4\zeta_{48}+a_5x_5\zeta_{44}
+a_6x_6\zeta_{40}+a_7x_7\zeta_{37}\\&
&+a_8x_8\zeta_{33}+a_9x_9\zeta_{30}+a_{10}x_{10}\zeta_{28}+a_{11}x_{11}\zeta_{25}+a_{12}x_{12}\zeta_{23}
+a_{13}x_{13}\zeta_{21}+a_{14}x_{14}\zeta_{19}\\&
&+a_{15}x_{15}\zeta_{17}+a_{16}x_{16}\zeta_{16}+a_{18}x_{18}\zeta_{14}+a_{19}x_{19}\zeta_{12}
+a_{20}x_{20}\zeta_{11}+a_{21}x_{21}\zeta_{10}
+a_{23}x_{23}\zeta_8\\& &
+a_{24}x_{24}\zeta_7+a_{26}x_{26}\zeta_6+a_{27}x_{27}\zeta_5
+a_{29}x_{29}\zeta_4+a_{32}x_{32}\zeta_3+a_{35}x_{35}\zeta_2+a_{40}x_{40}\zeta_1,\hspace{0.4cm}(3.115)\end{eqnarray*}
where $a_i\in\mbb{C}$. Then (2.52) and (3.92) yield
\begin{eqnarray*}0&=&E_{\al_1}(\vt)=-(a_6+a_8)x_6\zeta_{33}-(a_9+a_{11})x_9\zeta_{25}
-(a_{10}+a_{13})x_{10}\zeta_{21}\\ &
&-(a_{12}+a_{16})x_{12}\zeta_{16}
-(a_{14}+a_{19})x_{14}\zeta_{12}+(a_{40}-a_{35})x_{35}\zeta_1,\hspace{3.2cm}(3.116)
\end{eqnarray*}
equivalently,
$$a_8=-a_6,\;\;a_{11}=-a_9,\;\;a_{13}=-a_{10},\;\;a_{16}=-a_{12},\;\;a_{19}=-a_{14},\;\;a_{40}=a_{35}.\eqno(3.117)$$
Moreover, (2.53) and (3.93) give
\begin{eqnarray*}0&=&E_{\al_2}(\vt)=(a_5+a_7)x_5\zeta_{37}+(a_6+a_9)x_6\zeta_{30}
+(a_8+a_{11})x_8\zeta_{25}\\ & &+(a_{20}-a_{23})x_{20}\zeta_{16}
+(a_{24}-a_{26})x_{24}\zeta_6+(a_{27}-a_{29})x_{27}\zeta_4,\hspace{3.4cm}(3.118)
\end{eqnarray*}
equivalently,
$$a_7=-a_5,\;\;a_9=-a_6,\;\;a_{11}=-a_8,\;\;a_{23}=a_{20},\;\;a_{26}=a_{24},\;\;a_{29}=a_{27}.\eqno(3.119)$$
Furthermore, (2.54) and (3.94) imply
\begin{eqnarray*}0&=&E_{\al_3}(\vt)=-(a_5+a_6)x_5\zeta_{40}-(a_7+a_9)x_7\zeta_{30}
-(a_{13}+a_{15})x_{13}\zeta_{17}\\ &
&-(a_{16}+a_{18})x_{16}\zeta_{14}
-(a_{19}+a_{21})x_{19}\zeta_{10}+(a_{35}-a_{32})x_{32}\zeta_2,\hspace{3.2cm}(3.120)
\end{eqnarray*}
equivalently,
$$a_6=-a_5,\;\;a_9=-a_7,\;\;a_{15}=-a_{13},\;\;a_{18}=-a_{16},\;\;a_{21}=-a_{19},\;\;a_{35}=a_{32}.\eqno(3.121)$$
Note that (2.55) and (3.95) yield
\begin{eqnarray*}0&=&E_{\al_4}(\vt)=(a_4+a_5)x_4\zeta_{44}+(a_9-a_{10})x_9\zeta_{28}
\\ &
&+(a_{11}-a_{13})x_{11}\zeta_{21}+(a_{18}-a_{20})x_{18}\zeta_{11}
+(a_{29}-a_{32})x_{29}\zeta_3,\hspace{3.3cm}(3.122)
\end{eqnarray*}
equivalently,
$$a_5=-a_4,\;\;a_{10}=a_9,\;\;a_{13}=a_{11},\;\;a_{20}=a_{18},\;\;a_{32}=a_{29}.\eqno(3.123)$$
In addition, (2.56) and (3.96) say that
\begin{eqnarray*}0&=&E_{\al_5}(\vt)=(a_3+a_4)x_3\zeta_{48}+(a_{10}-a_{12})x_{10}\zeta_{23}
+(a_{13}-a_{16})x_{13}\zeta_{16}\\ &
&+(a_{15}-a_{18})x_{15}\zeta_{14}+(a_{24}-a_{27})x_{24}\zeta_5+(a_{26}-a_{29})x_{26}\zeta_4,\hspace{3.3cm}(3.124)
\end{eqnarray*}
equivalently,
$$a_4=-a_3,\;\;a_{12}=a_{10},\;\;a_{16}=a_{13},\;\;a_{18}=a_{15},\;\;a_{27}=a_{24},\;\;a_{29}=a_{26}.\eqno(3.125)$$
Besides, (2.57) and (3.97) tell us  that
\begin{eqnarray*}0&=&E_{\al_6}(\vt)=(a_2+a_3)x_2\zeta_{53}+(a_{12}-a_{14})x_{12}\zeta_{19}
+(a_{16}-a_{19})x_{16}\zeta_{12}\\ &
&+(a_{18}-a_{21})x_{18}\zeta_{10}+(a_{20}-a_{24})x_{20}\zeta_7+(a_{23}-a_{26})x_{23}\zeta_6,\hspace{3.4cm}(3.126)
\end{eqnarray*}
equivalently,
$$a_3=-a_2,\;\;a_{14}=a_{12},\;\;a_{19}=a_{16},\;\;a_{21}=a_{18},\;\;a_{24}=a_{20},\;\;a_{26}=a_{23}.\eqno(3.127)$$
Finally, (2.58), (3.98) and (3.114) show  that
$$0=E_{\al_7}(\vt)=(2a_1+a_2)x_1\zeta_{59}\lra a_2=-2a_1.\eqno(3.128)$$
Thus we have:
\begin{eqnarray*}\hspace{1cm}&&
2a_1=-a_2=a_3=-a_4=a_5=-a_6=-a_7=a_8=a_9=a_{10}=-a_{11}\\
&&=a_{12}=-a_{13}=a_{14}=a_{15}=-a_{16}=a_{18}=-a_{19}=a_{20}=a_{21}\\
&&=a_{23}=a_{24}=a_{26}=a_{27}=a_{29}
=a_{32}=a_{35}=a_{40}.\hspace{4.3cm}(3.129)\end{eqnarray*} Hence we
have the singular vector
\begin{eqnarray*}\vt&=&\frac{x_1}{2}\zeta_{64}-x_2\zeta_{59}+x_3\zeta_{53}-
x_4\zeta_{48}+x_5\zeta_{44}-x_6\zeta_{40}-x_7\zeta_{37}\\&
&+x_8\zeta_{33}+x_9\zeta_{30}+x_{10}\zeta_{28}-x_{11}\zeta_{25}+x_{12}\zeta_{23}
-x_{13}\zeta_{21}+x_{14}\zeta_{19}\\&
&+x_{15}\zeta_{17}-x_{16}\zeta_{16}+x_{18}\zeta_{14}-x_{19}\zeta_{12}
+x_{20}\zeta_{11}+x_{21}\zeta_{10} +x_{23}\zeta_8\\& &
+x_{24}\zeta_7+x_{26}\zeta_6+x_{27}\zeta_5
+x_{29}\zeta_4+x_{32}\zeta_3+x_{35}\zeta_2+x_{40}\zeta_1\hspace{3.8cm}(3.130)\end{eqnarray*}
of weight $\lmd_7$.

By Table 2, we want to find the singular vector of weight $\lmd_6$:
$$\vs=b_1\zeta_1\zeta_{19}+b_2\zeta_2\zeta_{12}+b_3\zeta_3\zeta_{10}+b_4\zeta_4\zeta_7+b_5\zeta_5\zeta_6,\qquad
b_j\in\mbb{C}.\eqno(3.131)$$ According (3.92),
$$0=E_{\al_1}(\vs)=-(b_1+b_2)\zeta_1\zeta_{12}\lra
b_2=-b_1.\eqno(3.132)$$ Moreover, (3.93) implies
$$0=E_{\al_2}(\vs)=(b_4+b_5)\zeta_4\zeta_6\lra
b_5=-b_4.\eqno(3.133)$$ Furthermore, (3.94) gives
$$0=E_{\al_3}(\vs)=-(b_2+b_3)\zeta_2\zeta_{10}\lra
b_3=-b_2.\eqno(3.134)$$ In addition, (3.95) yields
$$0=E_{\al_4}(\vs)=(b_3+b_4)\zeta_3\zeta_7\lra
b_4=-b_3.\eqno(3.135)$$ Observe that (3.133) implies
$E_{\al_5}(\vs)=0$ by (3.96) and $E_{\al_6}(\vs)=E_{\al_7}(\vs)=0$
naturally hold by  (3.97), (3.98) and (3.131).  Thus we have the
singular vector
$$\vs=\zeta_1\zeta_{19}-\zeta_2\zeta_{12}+\zeta_3\zeta_{10}-\zeta_4\zeta_7+\zeta_5\zeta_6\eqno(3.136)$$
of weight $\lmd_6$.

Observe that we have the module isomorphism $\sgm:{\cal G}^{E_7}\rta
W$ is determined by
$$\sgm(E_{\al_1})=-\zeta_{58},\;\;\sgm(E_{\al_2})=-\zeta_{57},\;\;\sgm(E_{\al_3})=\zeta_{63},
\;\;\sgm(E_{\al_4})=\zeta_{62},\eqno(3.137)$$
$$\sgm(E_{\al_5})=-\zeta_{61},\;\;\sgm(E_{\al_6})=\zeta_{60},\;\;\sgm(E_{\al_7})=-\zeta_{59},
\;\;\sgm(\al_1)=-\zeta_{70},\eqno(3.138)$$
$$\sgm(\al_2)=\zeta_{68},\;\;\sgm(\al_3)=\zeta_{69},\;\;\sgm(\al_4)=-\zeta_{67},\;\;
\sgm(\al_5)=\zeta_{66},\;\;\sgm(\al_6)=-\zeta_{65},\eqno(3.139)$$
$$\sgm(\al_7)=\zeta_{64},\;\;\sgm(E_{-\al_1})=\zeta_{76},\;\;
\sgm(E_{-\al_2})=-\zeta_{77},\;\;\sgm(E_{-\al_3})=-\zeta_{71},\eqno(3.140)$$
$$\sgm(E_{-\al_4})=\zeta_{72},\;\;\sgm(E_{-\al_5})=-\zeta_{73},\;\;\sgm(E_{-\al_6})=\zeta_{74},
\;\;\sgm(E_{-\al_7})=-\zeta_{75}.\eqno(3.141)$$

Next we use the above isomorphism to find the invariant of the form:
$$\eta=\sum_{i=1}^{63}c_i\zeta_i\zeta_{134-i}+\sum_{64\leq r\leq s\leq
70}c_{r,s}\zeta_r\zeta_s\eqno(3.142)$$ with $c_i,c_{r,s}\in\mbb{C}$.
 According to (3.92),
(3.107) and (3.108),
\begin{eqnarray*}0&=&E_{\al_1}(\eta)=(c_1-c_2)\zeta_1\zeta_{132}+(c_{12}-c_{19})\zeta_{12}\zeta_{115}
+(c_{16}-c_{23})\zeta_{16}\zeta_{111}+(c_{18}-c_{26})\zeta_{18}\zeta_{108}
\\ & &+(c_{21}-c_{28})\zeta_{21}\zeta_{106}+(c_{25}-c_{30})\zeta_{25}\zeta_{104}
+(c_{24}-c_{31})\zeta_{24}\zeta_{103}+(c_{29}-c_{34})\zeta_{29}\zeta_{100}\\
&&+(c_{27}-c_{35})\zeta_{27}\zeta_{99}+(c_{32}-c_{39})\zeta_{32}\zeta_{95}
+(c_{33}-c_{40})\zeta_{33}\zeta_{94}
+(c_{36}-c_{43})\zeta_{36}\zeta_{91}\\
&&+(c_{38}-c_{45})\zeta_{38}\zeta_{89}+(c_{27}-c_{35})\zeta_{27}\zeta_{99}+(c_{33}-c_{40})\zeta_{33}\zeta_{94}
+(c_{36}-c_{43})\zeta_{36}\zeta_{91}\\
&&+(c_{38}-c_{45})\zeta_{38}\zeta_{89}+(c_{42}-c_{50})\zeta_{42}\zeta_{84}+(c_{47}-c_{56})\zeta_{47}\zeta_{78}
+(c_{52}-c_{63})\zeta_{52}\zeta_{71}\\
& &-\zeta_{58}[(4c_{70,70}+c_{69,70}-c_{58})\zeta_{70}+\sum_{
i\in\ol{64,69}}(2c_{i,70}+c_{i,69})\zeta_i+c_{69,69}\zeta_{69}].\hspace{2cm}(3.143)\end{eqnarray*}
So
$$c_2=c_1,\;\;c_{19}=c_{12},\;\;c_{23}=c_{16},\;\;c_{26}=c_{18},\;\;c_{28}=c_{21},\;\;c_{30}=c_{25},
\eqno(3.144)$$
$$c_{31}=c_{24},\;\;c_{34}=c_{29},\;\;c_{35}=c_{27},\;\;c_{39}=c_{32},\;\;c_{40}=c_{33},\eqno(3.145)$$
$$c_{43}=c_{36},\;\;c_{45}=c_{38},
\;\;c_{50}=c_{42},\;\;c_{56}=c_{47},\;\;c_{63}=c_{52}.\eqno(3.146)$$

Similarly, (3.93), (3.107) and (3.109) imply that the equation
$E_{\al_2}(\eta)=0$ gives
$$c_5=-c_4,\;\;c_7=-c_6,\;\;c_{11}=-c_8,\;\;c_{13}=-c_9,\;\;c_{33}=-c_{25},\;\;c_{38}=-c_{29},
\eqno(3.147)$$
$$c_{40}=-c_{30},\;\;c_{42}=-c_{32},\;\;c_{44}=-c_{37},\;\;c_{45}=-c_{34},\;\;c_{47}=-c_{36},
\eqno(3.148)$$
$$c_{49}=-c_{41},\;\;c_{50}=-c_{39},
\;\;c_{55}=-c_{46},\;\;c_{56}=-c_{43},\;\;c_{62}=-c_{51}.\eqno(3.149)$$
 Moreover, (3.94), (3.107), (3.110) and the equation $E_{\al_3}(\eta)=0$
 yield
$$c_3=c_2,\;\;c_{12}=c_{10},\;\;c_{16}=c_{14},\;\;c_{18}=c_{15},\;\;c_{21}=c_{17},
\;\;c_{27}=c_{22},\eqno(3.150)$$
$$c_{24}=c_{20},\;\;c_{37}=c_{30},\;\;c_{41}=c_{34},\;\;c_{44}=c_{40},\;\;c_{46}=c_{39},\eqno(3.151)$$
$$c_{49}=c_{45},\;\;c_{51}=c_{43},
\;\;c_{55}=c_{50},\;\;c_{58}=c_{52},\;\;c_{62}=c_{56}.
\eqno(3.152)$$

Note that (3.95), (3.107), (3.111) and the equation
$E_{\al_4}(\eta)=0$  imply
$$c_4=-c_3,\;\;c_{10}=-c_7,\;\;c_{14}=-c_{11},\;\;c_{15}=-c_{13},\;\;c_{25}=-c_{21},\;\;c_{29}=-c_{24},
\eqno(3.153)$$
$$c_{30}=-c_{28},\;\;c_{32}=-c_{27},\;\;c_{34}=-c_{31},\;\;c_{39}=-c_{35},\;\;c_{48}=-c_{44},
\eqno(3.154)$$
$$c_{52}=-c_{47},\;\;c_{54}=-c_{49},
\;\;c_{57}=-c_{51},\;\;c_{61}=-c_{55},\;\;c_{63}=-c_{56}.\eqno(3.155)$$
Furthermore, (3.96), (3.107) and (3.112) give
$$c_6=-c_4,\;\;c_7=-c_5,\;\;c_{17}=-c_{14},\;\;c_{20}=-c_{15},\;\;c_{21}=-c_{16},\;\;c_{24}=-c_{18},
\eqno(3.156)$$
$$c_{28}=-c_{23},\;\;c_{31}=-c_{26},\;\;c_{36}=-c_{32},\;\;c_{43}=-c_{39},\;\;c_{47}=-c_{42},
\eqno(3.157)$$
$$c_{51}=-c_{46},\;\;c_{53}=-c_{48},\;\;c_{56}=-c_{50},\;\;c_{60}=-c_{54},\;\;c_{62}=-c_{55}.\eqno(3.158)$$

In addition, (3.97), (3.107) and (3.113) show
$$c_8=-c_6,\;\;c_{11}=-c_7,\;\;c_{14}=-c_{10},\;\;c_{16}=-c_{12},\;\;c_{22}=-c_{20},\;\;c_{23}=-c_{19},
\eqno(3.159)$$
$$c_{27}=-c_{24},\;\;c_{32}=-c_{29},\;\;c_{35}=-c_{31},\;\;c_{39}=-c_{34},\;\;c_{42}=-c_{38},
\eqno(3.160)$$
$$c_{46}=-c_{41},\;\;c_{50}=-c_{45},\;\;c_{55}=-c_{49},\;\;c_{59}=-c_{53},\;\;c_{61}=-c_{54}.\eqno(3.161)$$
Finally, (3.98), (3.107) and (3.114) yield
$$c_9=-c_8,\;\;c_{13}=-c_{11},\;\;c_{15}=-c_{14},\;\;c_{18}=-c_{16},\;\;c_{20}=-c_{17},\;\;c_{24}=-c_{21},
\eqno(3.162)$$
$$c_{26}=-c_{23},\;\;c_{29}=-c_{25},\;\;c_{31}=-c_{28},\;\;c_{34}=-c_{30},\;\;c_{38}=-c_{33},
\eqno(3.163)$$
$$c_{41}=-c_{37},\;\;c_{45}=-c_{40},\;\;c_{49}=-c_{44},\;\;c_{54}=-c_{48},\;\;c_{60}=-c_{53}.
 \eqno(3.164)$$

First we have the following solution of (3.144)-(3.164):
\begin{eqnarray*}&&
c_1=c_2=c_3=-c_4=c_5=c_6=-c_7=-c_8=c_9=c_{10}=c_{11}=c_{12}\\
&&=-c_{13}=-c_{14}=c_{15}=-c_{16}=c_{17}=c_{18}=c_{19}=-c_{20}=c_{21}=c_{22}=-c_{23}\\
&&=-c_{24}=-c_{25}=c_{26}=c_{27}=c_{28}=c_{29}=-c_{30}=-c_{31}=-c_{32}=c_{33}\\
&&=c_{34}=c_{35}=c_{36}=-c_{37}=-c_{38}=-c_{39}=c_{40}=c_{41}=c_{42}=c_{43}\\
&&=c_{44}=-c_{45}=-c_{46}=-c_{47}=-c_{48}=-c_{49}=c_{50}=c_{51}=c_{52}=c_{53}\\
&&=c_{54}=c_{55}=-c_{56}=-c_{57}=c_{58}=-c_{59}=-c_{60}=-c_{61}=-c_{62}=c_{63}.\hspace{1.5cm}(3.165)
\end{eqnarray*}
We want to find $c_{r,s}$ by (3.137)-(3.141). Define a bilinear form
$\la\cdot,\cdot\ra$ on ${\cal G}^{E_7}$ by
$$\la E_\al,E_\be\ra=-\dlt_{\al+\be,0}
\qquad\for\;\;\al,\be\in\Phi_{E_7},\eqno(3.166)$$
$$\la h,h'\ra=(h,h')\qquad\for\;\;h,h'\in H_{E_7}.\eqno(3.167)$$
Then $\la\cdot,\cdot\ra$ is an invariant bilinear form (e.g., cf.
[22], [42]). Its corresponding quadratic invariant polynomial over
${\cal G}^{E_7}$ is:
$$\Omega=\sum_{i=1}^7\al_i\lmd_i-2\sum_{\al\in\Phi_{E_7}^+}E_\al
E_{-\al},\eqno(3.168)$$ where the fundamental weights:
$$\lmd_1=2\al_1+2\al_2+3\al_3+4\al_4+3\al_5+2\al_6+\al_7,\eqno(3.169)$$
$$\lmd_2=\frac{1}{2}(4\al_1+7\al_2+8\al_3+12\al_4+9\al_5+6\al_6+3\al_7),
\eqno(3.170)$$
$$\lmd_3=3\al_1+4\al_2+6\al_3+8\al_4+6\al_5+4\al_6+2\al_7,\eqno(3.171)$$
$$\lmd_4=4\al_1+6\al_2+8\al_3+12\al_4+9\al_5+6\al_6+3\al_7,\eqno(3.172)$$
$$\lmd_5=\frac{1}{2}(6\al_1+9\al_2+12\al_3+18\al_4+15\al_5+10\al_6+5\al_7),
\eqno(3.173)$$
$$\lmd_6=2\al_1+3\al_2+4\al_3+6\al_4+5\al_5+4\al_6+2\al_7,\eqno(3.174)$$
$$\lmd_7=\frac{1}{2}(2\al_1+3\al_2+4\al_3+6\al_4+5\al_5+4\al_6+3\al_7).
\eqno(3.175)$$ In particular,
\begin{eqnarray*}\sum_{i=1}^7\al_i\lmd_i&=&\al_1(2\al_1+4\al_2+6\al_3+8\al_4+6\al_5+4\al_6+2\al_7)+
\al_2(7\al_2/2+8\al_3+12\al_4\\ &
&+9\al_5+6\al_6+3\al_7)+\al_3(6\al_3+16\al_4+12\al_5+8\al_6+4\al_7)+\al_4(12\al_4+18\al_5\\
& &+12\al_6+6\al_7)+\al_5(15\al_5/2+10\al_6+5\al_7)
+4\al_6(\al_6+\al_7)+3\al_7^2/2.\hspace{0.3cm}(3.176)\end{eqnarray*}
Note that (3.137)-(3.141) say that
$$\sgm(\sum_{i=1}^7E_{\al_i}E_{-\al_i})=\zeta_{57}\zeta_{77}-\zeta_{58}\zeta_{76}-\zeta_{63}\zeta_{71}+\sum_{i=59}^{62}
\zeta_i\zeta_{134-i},\eqno(3.177)$$
\begin{eqnarray*}\sgm(2\sum_{i=1}^7\al_i\lmd_i)&=&
4\zeta_{70}(\zeta_{70}-2\zeta_{68}-3\zeta_{69}+4\zeta_{67}-3\zeta_{66}+2\zeta_{65}-\zeta_{64})\\
& &+
\zeta_{68}(7\zeta_{68}+16\zeta_{69}-24\zeta_{67}+18\zeta_{66}-12\zeta_{65}+6\zeta_{64})+4\zeta_{69}(3\zeta_{69}\\
&&-8\zeta_{67}+6\zeta_{66}
-4\zeta_{65}+2\zeta_{64})+4\zeta_{67}(6\zeta_{67}-9\zeta_{66}+6\zeta_{65}-3\zeta_{64})\\
& &+\zeta_{66}(15\zeta_{66}-20\zeta_{65}+10\zeta_{64})
+8\zeta_{65}(\zeta_{65}-\zeta_{64})+3\zeta_{64}^2.\hspace{2.2cm}(3.178)\end{eqnarray*}
Denote
$$I=\{\ol{1,6},\ol{9,12},15,\ol{17,19},21,22,\ol{26.29},\ol{33,36},\ol{40,44},\ol{50,55},58,63\}.\eqno(3.179)$$
Taking $c_1=4$, we obtain the following explicit formula of quartic
invariant:\begin{eqnarray*}\eta &=&4(\sum_{i\in
I}\zeta_i\zeta_{134-i}-\sum_{r\in\ol{7,63}\setminus
I}\zeta_r\zeta_{134-r})+
4\zeta_{70}(\zeta_{70}-2\zeta_{68}-3\zeta_{69}+4\zeta_{67}-3\zeta_{66}\\
& &+2\zeta_{65}-\zeta_{64})+
\zeta_{68}(7\zeta_{68}+16\zeta_{69}-24\zeta_{67}+18\zeta_{66}-12\zeta_{65}+6\zeta_{64})\\
&&+4\zeta_{69}(3\zeta_{69}-8\zeta_{67}+6\zeta_{66}
-4\zeta_{65}+2\zeta_{64})+4\zeta_{67}(6\zeta_{67}-9\zeta_{66}+6\zeta_{65}\\
& &-3\zeta_{64})+\zeta_{66}(15\zeta_{66}-20\zeta_{65}+10\zeta_{64})
+8\zeta_{65}(\zeta_{65}-\zeta_{64})+3\zeta_{64}^2.\hspace{2.8cm}(3.180)\end{eqnarray*}

\section{Proof of Main Theorem}

In this section, we want to prove the main theorem of this paper.
\psp

{\bf Theorem 4.1}. {\it Any singular vector in ${\cal A}$ is a
polynomial in $x_1,\;\zeta_1,\;\vt,\;\vs$ and $\eta$. Let
$L(n_1,n_2,n_3,n_4,n_5)$ be the irreducible submodule generated by
$\zeta_1^{n_1}\vs^{n_2}x_1^{n_3}\vt^{n_4}\eta^{n_5}$ with highest
weight $n_1\lmd_1+n_2\lmd_6+(n_3+n_4)\lmd_7$. Then
$${\cal A}=\bigoplus_{n_1,n_2,n_3,n_4,n_5=0}^\infty
L(n_1,n_2,n_3,n_4,n_5).\eqno(4.1)$$ In particular,
\begin{eqnarray*}& &(1-q)^{55}\sum_{n_1,n_2,n_3,n_4=0}^\infty(\mbox{dim}\:
V(n_1\lmd_1+n_2\lmd_6+(n_3+n_4)\lmd_7))q^{2n_1+4n_2+n_3+3n_4}\\
&=&1+q+q^2+q^3.\hspace{11.3cm}(4.2)\end{eqnarray*}
  Let ${\cal D}$ be the invariant differential operator obtained
  from $\eta$ by changing $x_i$ to $\ptl_{x_i}$. Then
  $$\sum_{n_1,n_2,n_3=0}^\infty(L(n_1,n_2,n_3,0,0)+L(n_1,n_2,n_3,1,0))\subset
  \{f\in{\cal A}\mid {\cal D}(f)=0\}.\eqno(4.3)$$}

{\it Proof}. Let $f$ be a singular vector in ${\cal A}$. To
eliminate the extra variables from $f$, we need a certain change of
variables. According to (3.7), (3.18)-(3.24), (3.26)-(3.28), (3.30),
(3.32), (3.33), (3.37), (3.39), (3.41), (3.44), (3.46), (3.49),
(3.53), (3.56), (3.60), (3.64), (3.69), (3.75), and (3.136), we have
$$x_1x_{17}=\zeta_1-x_2x_{14}-x_3x_{12}-x_4x_{10}-x_5x_9+x_6x_7.
\eqno(4.4)$$
$$x_1x_{22}=\zeta_2-x_2x_{19}-x_3x_{16}-x_4x_{13}-x_5x_{11}+x_7x_8,
\eqno(4.5)$$
$$x_1x_{25}=\zeta_3-x_2x_{21}-x_3x_{18}-x_4x_{15}-x_6x_{11}+x_8x_9,
\eqno(4.6)$$
$$x_1x_{28}=-\zeta_4-x_2x_{24}-x_3x_{20}+x_5x_{15}-x_6x_{13}+x_8x_{10},
\eqno(4.7)$$
$$x_1x_{30}=\zeta_5-x_2x_{26}-x_3x_{23}-x_7x_{15}+x_9x_{13}-x_{10}x_{11},
\eqno(4.8)$$
$$x_1x_{31}=\zeta_6-x_2x_{27}+x_4x_{20}+x_5x_{18}-x_6x_{16}+x_8x_{12},
\eqno(4.9)$$
$$x_1x_{33}=-\zeta_7-x_2x_{29}+x_4x_{23}-x_7x_{18}+x_9x_{16}-x_{11}x_{12},
\eqno(4.10)$$
$$x_1x_{34}=-\zeta_8+x_3x_{27}+x_4x_{24}+x_5x_{21}-x_6x_{19}+x_8x_{14},
\eqno(4.11)$$
$$x_1x_{36}=\zeta_{10}-x_2x_{32}-x_5x_{23}-x_7x_{20}+x_{10}x_{16}-x_{12}x_{13},
\eqno(4.12)$$
$$x_1x_{37}=\zeta_{11}+x_3x_{29}+x_4x_{26}-x_7x_{21}+x_9x_{19}-x_{11}x_{14},
\eqno(4.13)$$
$$x_1x_{38}=-\zeta_{12}-x_2x_{35}+x_6x_{23}+x_9x_{20}-x_{10}x_{18}+x_{12}x_{15},
\eqno(4.14)$$
$$x_1x_{39}=\zeta_{14}+x_3x_{32}-x_5x_{26}-x_7x_{24}+x_{10}x_{19}-x_{13}x_{14},
\eqno(4.15)$$
$$x_1x_{41}=\zeta_{16}+x_3x_{35}+x_6x_{26}+x_9x_{24}-x_{10}x_{21}+x_{14}x_{15},
\eqno(4.16)$$
$$x_1x_{42}=\zeta_{17}-x_4x_{32}-x_5x_{29}-x_7x_{27}+x_{12}x_{19}-x_{14}x_{16},
\eqno(4.17)$$
\begin{eqnarray*}\hspace{1cm}x_1\zeta_1x_{43}&=&
-\zeta_1(x_1x_{43}+x_2x_{40}+x_8x_{23}+x_{11}x_{20}-x_{13}x_{18}+x_{15}x_{16})
\\ & &+\vs+\zeta_2\zeta_{12}-\zeta_3\zeta_{10}+\zeta_4\zeta_7-\zeta_5\zeta_6,\hspace{5.7cm}(4.18)\end{eqnarray*}
$$x_1x_{44}=-\zeta_{21}-x_4x_{35}+x_6x_{29}+x_9x_{27}-x_{12}x_{21}+x_{14}x_{18},
\eqno(4.19)$$
$$x_1x_{45}=-\zeta_{23}+x_3x_{40}-x_8x_{26}-x_{11}x_{24}+x_{13}x_{21}-x_{15}x_{19},
\eqno(4.20)$$
$$x_1x_{46}=\zeta_{25}+x_5x_{35}+x_6x_{32}+x_{10}x_{27}-x_{12}x_{24}+x_{14}x_{20},
\eqno(4.21)$$
$$x_1x_{47}=\zeta_{28}-x_4x_{40}-x_8x_{29}-x_{11}x_{27}+x_{16}x_{21}-x_{18}x_{19},
\eqno(4.22)$$
$$x_1x_{48}=-\zeta_{30}+x_5x_{40}-x_8x_{32}-x_{13}x_{27}+x_{16}x_{24}-x_{19}x_{20},
\eqno(4.23)$$
$$x_1x_{49}=-\zeta_{33}+x_7x_{35}+x_9x_{32}-x_{10}x_{29}+x_{12}x_{26}-x_{14}x_{23},
\eqno(4.24)$$
$$x_1x_{50}=\zeta_{37}-x_6x_{40}-x_8x_{35}+x_{15}x_{27}-x_{18}x_{24}+x_{20}x_{21},
\eqno(4.25)$$
$$x_1x_{51}=-\zeta_{40}+x_7x_{40}-x_{11}x_{32}+x_{13}x_{29}-x_{16}x_{26}+x_{19}x_{23},
\eqno(4.26)$$
$$x_1x_{52}=\zeta_{44}-x_9x_{40}-x_{11}x_{35}-x_{15}x_{29}+x_{18}x_{26}-x_{21}x_{23},
\eqno(4.27)$$
$$x_1x_{53}=\zeta_{48}+x_{10}x_{40}+x_{13}x_{35}+x_{15}x_{32}-x_{20}x_{26}+x_{23}x_{24},
\eqno(4.28)$$
$$x_1x_{54}=\zeta_{53}-x_{12}x_{40}-x_{16}x_{35}-x_{18}x_{32}+x_{20}x_{29}-x_{23}x_{27},
\eqno(4.29)$$
$$x_1x_{55}=\zeta_{59}+x_{14}x_{40}+x_{19}x_{35}+x_{21}x_{32}-x_{24}x_{29}+x_{26}x_{27}.
\eqno(4.30)$$

Moreover, (3.7) and (3.18)-(3.87) imply that when $x_i=0$ for
$i\in\ol{3,54}$,
$$\vt=\frac{x_1}{2}(x_1x_{56}-x_2x_{55}),\eqno(4.31)$$
which is the homogeneous part of $\vt$ with degree 1 in $x_{55}$ and
$x_{56}$ (cf. (3.130)), and
$$\eta=3x_1^2x_{56}^2-6x_1x_2x_{55}x_{56}-5x_2^2x_{55}^2,\eqno(4.32)$$
which is  the homogeneous part of $\eta$ with degree 2 in $x_{55}$
and $x_{56}$ (cf. (3.180)). Under the assumption, we substitute
$x_1x_{56}=2x_1^{-1}\vt+x_2x_{55}$ into (4.32) and obtain
$$\eta=-8x_2^2x_{55}^2+12x_1^{-2}\vt^2.\eqno(4.33)$$
This shows that $f$ can be written as a function in
$\{x_i,\vt,\eta\mid i\in\ol{1,54}\}$. Moreover, using (4.4)-(4.30),
we obtain that $f$ is a function $f_1$ in
\begin{eqnarray*}\hspace{1cm}& &\{x_r,\zeta_s,\vt,\vs,\eta\mid r\in \{\ol{1,35},40\}\setminus
\{17,22,25,28,30,31,33,34\},\;s\in\{\ol{1,17},\\
&&21,23,25,28,30,33,37,40,44,48,53\}\setminus\{9,13,15\}\}.\hspace{4.1cm}(4.34)\end{eqnarray*}

Set
$$H_{E_6}=\sum_{i=1}^6\mbb{C}\al_i.\eqno(4.35)$$
Then
$$\Phi_{E_6}=\Phi_{E_7}\bigcap H_{E_6}\eqno(4.36)$$
is a root system of $E_6$. Moreover,
$${\cal G}^{E_6}=H_{E_6}+\sum_{\al\in
\Phi_{E_6}}\mbb{C}E_\al\eqno(4.37)$$ forms a simple Lie algebra of
type $E_6$. Since ${\cal G}^{E_6}$ is a subalgebra of ${\cal
G}^{E_7}$ and $H_{E_6}\subset H_{E_7}$, $\zeta_1$ is also a
singular vector for ${\cal G}^{E_6}$ whose weight
$\lmd_1|_{H_{E_6}}$ is the first fundamental weight of ${\cal
G}^{E_6}$. Thus $\zeta_1$ generates a 27-dimensional irreducible
 ${\cal G}^{E_6}$-module $U$. By  (3.18)-(3.24), (3.26)-(3.28), (3.30),
(3.32), (3.33), (3.35), (3.37), (3.39), (3.41), (3.44), (3.46),
(3.49), (3.53), (3.56), (3.60), (3.64), (3.69), and (3.75), we
have
\begin{eqnarray*}\hspace{1cm}U&=&\sum_{i=1}^8\mbb{C}\zeta_i+\sum_{13,15\neq
r\in\ol{10,17}}\mbb{C}\zeta_r+\mbb{C}\zeta_{19}+\mbb{C}\zeta_{21}+\mbb{C}\zeta_{23}+\mbb{C}\zeta_{25}+\mbb{C}\zeta_{28}
+\mbb{C}\zeta_{30}\\
&&+\mbb{C}\zeta_{33}+\mbb{C}\zeta_{37}+\mbb{C}\zeta_{40}
+\mbb{C}\zeta_{44}+\mbb{C}\zeta_{48}
+\mbb{C}\zeta_{53}+\mbb{C}\zeta_{59}.\hspace{3.5cm}(4.38)\end{eqnarray*}
Recalling the fact that
$$[E_{\al_7},E_{-\al_i}]=0\qquad\for\;\;i\in\ol{1,6},\eqno(4.39)$$
we have $E_{\al_7}|_U=0$. Thus
$$E_\gm|_U=0\qquad\for\;\;\gm\in \Phi_{E_7}^+\setminus \Phi_{E_6}^+,\eqno(4.40)$$
because $\{E_{\al_1},...,E_{\al_7}\}$ generates all positive root
vectors of ${\cal G}^{E_7}$.

According (2.105)-(2.114) and (4.40), we have
$$0=E_{(1,1,2,3,2,1,1)}(f_1)=-x_1\ptl_{x_{20}}(f_1),\eqno(4.41)$$
$$0=E_{(1,1,2,2,2,2,1)}(f_1)=-x_1\ptl_{x_{21}}(f_1),\eqno(4.42)$$
$$0=E_{(1,2,2,3,2,1,1)}(f_1)=x_1\ptl_{x_{23}}(f_1),\eqno(4.43)$$
$$0=E_{(1,1,2,3,2,2,1)}(f_1)=x_1\ptl_{x_{24}}(f_1),\eqno(4.44)$$
$$0=E_{(1,2,2,3,2,2,1)}(f_1)=-x_1\ptl_{x_{26}}(f_1),\eqno(4.45)$$
$$0=E_{(1,1,2,3,3,2,1)}(f_1)=-x_1\ptl_{x_{27}}(f_1),\eqno(4.46)$$
$$0=E_{(1,2,2,3,3,2,1)}(f_1)=x_1\ptl_{x_{29}}(f_1),\eqno(4.47)$$
$$0=E_{(1,2,2,4,3,2,1)}(f_1)=-x_1\ptl_{x_{32}}(f_1),\eqno(4.48)$$
$$0=E_{(1,2,3,4,3,2,1)}(f_1)=-x_1\ptl_{x_{35}}(f_1),\eqno(4.49)$$
$$0=E_{(2,2,3,4,3,2,1)}(f_1)=-x_1\ptl_{x_{40}}(f_1).\eqno(4.50)$$
So $f_1$ is independent of
$\{x_{20},x_{21},x_{23},x_{24},x_{26},x_{27},x_{29},x_{32},x_{35},x_{40}\}$,
that is, $f_1$ is  a  function  in
\begin{eqnarray*}\hspace{2cm}& &\{x_r,\zeta_s,\vt,\vs\mid 17\neq r\in \ol{1,19};\;s\in\{\ol{1,17},21,23,\\
&&25,28,30,33,37,40,44,48,53\}\setminus\{9,13,15\}\}.\hspace{4.3cm}(4.51)\end{eqnarray*}

Next (2.58), (2.64), (2.70), (2.76), (2.81), (2.82), (2.85),
(2.87), (2.90), (2.92), (2.95), (2.96), (2.98)-(2.100), (2.102),
(2.103) and (4.40) imply
$$0=E_{(0,0,0,0,0,0,1)}(f_1)=x_1\ptl_{x_2}(f_1),\eqno(4.52)$$
$$0=E_{(0,0,0,0,0,1,1)}(f_1)=-x_1\ptl_{x_3}(f_1),\eqno(4.53)$$
$$0=E_{(0,0,0,0,1,1,1)}(f_1)=x_1\ptl_{x_4}(f_1),\eqno(4.54)$$
$$0=E_{(0,0,0,1,1,1,1)}(f_1)=-x_1\ptl_{x_5}(f_1),\eqno(4.55)$$
$$0=E_{(0,1,0,1,1,1,1)}(f_1)=x_1\ptl_{x_7}(f_1),\eqno(4.56)$$
$$0=E_{(0,0,1,1,1,1,1)}(f_1)=-x_1\ptl_{x_6}(f_1),\eqno(4.57)$$
$$0=E_{(1,0,1,1,1,1,1)}(f_1)=-x_1\ptl_{x_8}(f_1),\eqno(4.58)$$
$$0=E_{(0,1,1,1,1,1,1)}(f_1)=x_1\ptl_{x_9}(f_1),\eqno(4.59)$$
$$0=E_{(1,1,1,1,1,1,1)}(f_1)=x_1\ptl_{x_{11}}(f_1),\eqno(4.60)$$
$$0=E_{(0,1,1,2,1,1,1)}(f_1)=-x_1\ptl_{x_{10}}(f_1),\eqno(4.61)$$
$$0=E_{(1,1,1,2,1,1,1)}(f_1)=-x_1\ptl_{x_{13}}(f_1),\eqno(4.62)$$
$$0=E_{(0,1,1,2,2,1,1)}(f_1)=x_1\ptl_{x_{12}}(f_1),\eqno(4.63)$$
$$0=E_{(1,1,2,2,1,1,1)}(f_1)=-x_1\ptl_{x_{15}}(f_1),\eqno(4.64)$$
$$0=E_{(1,1,1,2,2,1,1)}(f_1)=x_1\ptl_{x_{16}}(f_1),\eqno(4.65)$$
$$0=E_{(0,1,1,2,2,2,1)}(f_1)=-x_1\ptl_{x_{14}}(f_1),\eqno(4.66)$$
$$0=E_{(1,1,2,2,2,1,1)}(f_1)=x_1\ptl_{x_{18}}(f_1),\eqno(4.67)$$
$$0=E_{(1,1,1,2,2,2,1)}(f_1)=-x_1\ptl_{x_{19}}(f_1).\eqno(4.68)$$
Hence $f_1$ is independent of $\{x_r,x_{18},x_{19}\mid\ol{2,16}\}$,
that is, $f_1$ is  a  function in
$$\{x_1,\zeta_s,\vt,\vs,\eta\mid
s\in\{\ol{1,17},21,23,25,28,30,33,37,40,44,48,53\}\setminus\{9,13,15\}\}.
\eqno(4.69)$$

Viewing $x_1,\vt,\vs,\eta$ as parameter constants, we can treat
$f_1$ as a ${\cal G}^{E_6}$ singular vector  on $U$, which does not
involve $\zeta_{19}$ and $\zeta_{59}$. According to our result in
[44], $f_1$ is a polynomial in $\zeta_1,\vs$ and a cubic invariant
involving $\zeta_{59}$. This shows that $f_1$ must be a function in
$\{x_1,\zeta_1,\vt,\vs,\eta\}$.

Note that $\vt$ and $\eta$ are algebraically independent over the
rational functions in $\{x_i\mid i\in\ol{3,54}\}$ by (4.31) and
(4.32). Moreover, (4.32) is an irreducible polynomial.  By (4.4) and
(4.18), $f=f_1$ must be a polynomial in
$\{x_1,\zeta_1,\vt,\vs,\eta\}$. Since each subspace of homogeneous
polynomials in ${\cal A}$ is a finite-dimensional ${\cal
G}^{E_7}$-module,  we get (4.1) by the Weyl's theorem of complete
reducibility and the fact that all finite-dimensional irreducible
${\cal G}^{E_7}$-modules are of highest-weight type. Moreover, (4.1)
gives
\begin{eqnarray*}& &(\sum_{r=0}^\infty q^4)\sum_{n_1,n_2,n_3,n_4=0}^\infty(\mbox{dim}\:
L(n_1,n_2,n_3,n_4,0))q^{2n_1+4n_2+n_3+3n_4}\\
&=&\frac{1}{(1-q)^{56}},\hspace{12cm}(4.70)\end{eqnarray*} that
is,
\begin{eqnarray*}& &\frac{1}{1-q^4}\sum_{n_1,n_2,n_3,n_4=0}^\infty(\mbox{dim}\:
L(n_1,n_2,n_3,n_4,0))q^{2n_1+4n_2+n_3+3n_4}\\
&=&\frac{1}{(1-q)^{56}},\hspace{11.9cm}(4.71)\end{eqnarray*} which
is equivalent to (4.2).

Observe that  ${\cal D}(\zeta_1^{m_1}\vs^{m_2}x_1^{m_3}\vt^{m_4})$
is also a singular vector of weight
$m_1\lmd_1+m_2\lmd_6+(m_3+m_4)\lmd_7$ if it is nonzero. So it is a
linear combination of the elements
$\zeta_1^{m_1}\vs^{m_2}x_1^{n_3}\vt^{n_4}\eta^{n_5}$ such that
$$n_3+n_4=m_3+m_4,\qquad n_3+3n_4+4n_5=m_3+3m_4-4.\eqno(4.72)$$
Thus
$$m_4=2+n_4+2n_5\geq 2.\eqno(4.73)$$
This shows
$${\cal D}(\zeta_1^{m_1}\vs^{m_2}x_1^{m_3})={\cal
D}(\zeta_1^{m_1}\vs^{m_2}x_1^{m_3}\vt)=0\eqno(4.74)$$ for any
nonnegative integers $m_1,m_2$ and $m_3$, equivalently, (4.3)
holds. The proof of our main theorem is completed. $\qquad\Box$

\bibliographystyle{amsplain}

\end{document}